\documentclass[11pt,pdftex]{article}
\pdfoutput = 1

\usepackage{psfrag}
\usepackage{color}
\usepackage{amsmath}
\usepackage{amsfonts}
\usepackage{graphicx}
\usepackage{amssymb}
\usepackage{rotating}
\usepackage{multirow}
\usepackage{subfigure}
\usepackage[normalem]{ulem} 
\usepackage{cancel} 

\usepackage{epstopdf}
\usepackage[unicode,colorlinks=true,linkcolor=blue,citecolor=blue,pagebackref=false,pdfborder={0 0 0},urlcolor=cyan]{hyperref}

\graphicspath{{FIGURES/}{figures/}}

\def\qed{\hbox{${\vcenter{\vbox{                        
   \hrule height 0.4pt\hbox{\vrule width 0.4pt height 6pt
   \kern5pt\vrule width 0.4pt}\hrule height 0.4pt}}}$}}
\newtheorem{theo}{Theorem}[section]
\newtheorem{prop}[theo]{Proposition}
\newtheorem{remark}{Remark}[section]
\newtheorem{lemm}[theo]{Lemma}
\newtheorem{coroll}[theo]{Corollary}
\newtheorem{assumpt}{Assumption}[section]

\newtheorem{proposition}{Proposition}[section]

\newcommand\bA{{\bf A}}
\newcommand\bp{{\bf p}}
\newcommand\bzeta{\boldsymbol \zeta}

\newcommand\bisect{{\rm bisect}}

\def\R{\mathbb{R}}

\def\lec{\lesssim}

\def\hsk{{\rm hSk}}
\def\vsk{{\rm vSk}}

\definecolor{blue}{rgb}{0,0,0.99}
\def\Bd{\color{black}}
\definecolor{dred}{rgb}{0.92,0,0}

\def\Bk{\color{black}}

\definecolor{dsky}{rgb}{0.1,1,1}

\newcommand\Aone{{\bf (A1)\ }}
\newcommand\Atwo{{\bf (A2)}}

\let\tilde\widetilde

\begin{document}


\title{BPX preconditioners for isogeometric analysis \\ using analysis-suitable T-splines}

\date{}
\author{Durkbin Cho\thanks{Department of Mathematics, Dongguk University, Pil-dong 3-ga, Jung-gu, Seoul, 100-715, South Korea. \newline \textit{E-mail address: } \texttt{durkbin@dongguk.edu}}
\and
Rafael V\'azquez\thanks{Istituto di Matematica Applicata e Tecnologie Informatiche `E. Magenes' del CNR, via Ferrata 1, 27100, Pavia, Italy. \textit{E-mail address: } \texttt{vazquez@imati.cnr.it} and
Institute of Mathematics, Ecole Polytechnique F\'ed\'erale de Lausanne, Station 8, 1015, Lausanne, Switzerland. \textit{E-mail address: } \texttt{rafael.vazquez@epfl.ch}
}
}



\maketitle

\begin{abstract}
  We propose and analyze optimal additive multilevel solvers for isogeometric discretizations of scalar elliptic problems for locally refined T-meshes. Applying the refinement strategy in \cite{Morgenstern_Peterseim} we can guarantee that the obtained T-meshes have a multilevel structure, and that the associated T-splines are analysis-suitable, for which we can define a dual basis and a stable projector. Taking advantage of the multilevel structure, we develop two BPX preconditioners: the first on the basis of local smoothing only for the functions affected by a newly added edge by bisection, and the second smoothing for all the functions affected after adding all the edges of the same level. We prove that both methods have optimal complexity, and present several numerical experiments to confirm our theoretical results, and also to compare the practical performance of the proposed preconditioners.
\end{abstract}



\section{Introduction}



The analysis and development of adaptive schemes is one of the most active areas of research in the context of isogeometric analysis (IGA), a recent methodology for the solution of partial differential equations with high continuity splines. The main idea of adaptive methods is to obtain a good accuracy of the solution with less computational effort by applying local mesh refinement, hence adaptive IGA schemes require the use of different spline spaces that break the tensor-product structure of B-splines. The most popular alternatives in the IGA research community are T-splines \cite{Sederberg_Zheng}, LR-splines \cite{LR-splines, Johannessen2013} and hierarchical splines \cite{Vuong_giannelli_juttler_simeon}.

In particular, in this paper we focus on T-splines, introduced by Sederberg et al. in \cite{Sederberg_Zheng} and applied in IGA for the first time in \cite{Bazilervs_Calo_Cottrell_Evans,Dorfel_Juttler_Simeon}. T-splines are constructed from a T-mesh, a rectangular tiling with hanging nodes, and T-spline blending functions are defined from their local knot vectors, which are computed from the tiling. The mathematical research on T-splines has been very active in recent years, and it has led to the introduction of analysis-suitable (or dual compatible) T-splines \cite{LZSHS12,BBCS12, BBSV2013,Li-Scott}, a subset of T-splines with good approximation properties that provide local linear independence.

A standard adaptive scheme based on mesh refinement can be written in a loop of the form \cite{CNX}
\begin{equation*}
\text{\bf SOLVE} \longrightarrow \text{\bf ESTIMATE} \longrightarrow \text{\bf MARK} \longrightarrow \text{\bf REFINE},
\end{equation*}
and suitable strategies for all the steps of the adaptive scheme are needed in order to guarantee its efficiency.
In this paper we focus on the solution of the linear system arising in the {\bf SOLVE} module for T-splines, and study the optimality of a suitable BPX preconditioner. Several domain decomposition preconditioners have been recently studied in the IGA context: overlapping Schwarz methods \cite{BePa13}, balancing domain decomposition by constraint (BDDC) methods \cite{BePa_BDDC,BPSWZ} and dual primal isogeometric tearing and interconnecting methods \cite{KPJ_FETIDP,PS_FETIDP}. Multilevel preconditioners for IGA have been extensively analyzed in the tensor-product setting: the BPX preconditioner in \cite{BHKS13}, and multigrid preconditioners \cite{Gahalaut_MG,Hofreither2015_1,Hofreither2015_2,Hofreither2016,Hofreither-Takacs,Manni_MG}. \Bd Moreover,  preconditioners based on fast solvers for Sylvester-like equations have been proposed in the recent paper \cite{Sangalli_Tani}.
\Bk

The T-splines obtained with the refinement strategy in \cite{Morgenstern_Peterseim} are analysis-suitable by construction, and present a multilevel structure, which makes them very appealing to apply multilevel preconditioners. In the present paper, we first present theoretically optimal multilevel preconditioners for IGA on {\it locally refined} T-meshes, extending the results in \cite{BHKS13} to T-splines.

For the study of the optimality of the BPX preconditioner we follow \cite{CNX}, writing the preconditioner in the framework of the parallel subspace correction (PSC) method. In this framework, the optimality follows from two basic properties: a stable space decomposition, and the strengthened Cauchy-Schwarz inequality. The proof of these two properties is the core of this work, and as a consequence we obtain that the BPX preconditioner gives a uniformly bounded condition number, which is independent of the mesh size $h$, but depends on the degree $p$ as in \cite{BHKS13}.


The construction of the BPX preconditioner as in \cite{CNX} is performed by adding a new edge to the T-mesh by bisection, and the new level is defined by the functions appearing or modified by the insertion of this edge. An alternative construction is also proposed, adding all the edges of the same generation at once, and defining the functions of the new level as those appearing or modified after the insertion of all edges. The theoretical optimality for this alternative construction, that we name \emph{macro} decomposition, follows from the previous one, but the numerical results show an improved performance.

The paper is organized as follows. In Section~2, we introduce the framework of parallel subspace correction (PSC) method and present its convergence theory based on the two properties mentioned above. We briefly review the basics of univariate/multivariate B-splines in Section~3. In Section~4, we give a new definition of T-meshes by bisections and discuss $\bp$-admissible T-meshes and their fundamental properties as in \cite{Morgenstern_Peterseim}. In Section~5, we construct a space decomposition on $\bp$-admissible T-meshes and then prove that the two aforementioned properties are satisfied. In Section~6 we obtain the optimality result for the BPX preconditioner, and also deal with the macro decomposition, showing that is also an efficient space decomposition for the purpose of implementation. Some numerical results that validate our theory are presented in Section~7.

\section{Preliminaries}
\subsection{Problem setting}
We are interested in the second order Laplacian with homogeneous Dirichlet boundary conditions,
\begin{equation}\label{model_prob}
-\Delta u = f \quad \mbox{\ in \ } \Omega, \qquad u = 0 \quad \mbox{\ on \ } \partial\Omega,
\end{equation}
where $\partial\Omega$ denotes the boundary of $\Omega$ and $f\in L^2(\Omega)$. The isogeometric approximation to the solution of \eqref{model_prob} is the function $u\in \mathcal{V}$ such that
\begin{equation*}
a(u,v) = ( f , v ) \quad \forall v\in \mathcal{V}
\end{equation*}
where
\[
a(u,v)=\int_\Omega \nabla u \cdot \nabla v~dx, \quad ( f, v ) = \int_\Omega f v~ dx,
\]
and ${\cal V}$ is the isogeometric discrete space.
Defining a linear operator $A:\mathcal{V}\rightarrow \mathcal{V}$ by
\[
(Au,v)=a(u,v), \ \forall u, v \in \mathcal{V}
\]
and also $b\in\mathcal{V}$ by $(b,v)=(f,v),\ \forall v\in \mathcal{V}$, we have to solve the linear operator equation
\begin{equation*}
Au=b
\end{equation*}
for some $u\in \mathcal{V}$.

In the rest of the paper, we will adopt the following compact notation. Given
two real numbers $a,b$ we write $a \lec b$,
when $a\leq C b$ for a generic constant $C$
independent of the knot vectors \Bd and the mesh size $h$ but depending on the spline degree $p$ and the geometric map ${\mathbf F}$ (defined below), \Bk and we write $a\approx b$ when $a \lec b$ and $b \lec a$.

\subsection{\Bd Preconditioned conjugate gradient method}
\Bd Let $B$ be a symmetric positive definite (SPD) operator. Applying it to both sides of $Au=b$, we get an equivalent equation
\[
BAu=Bb.
\]
The preconditioned conjugate gradient (PCG) method can be viewed as a conjugate gradient method applied to solving $BAu=Bb$ where $B$ is called a {\it preconditioner} (see, e.g., \cite{CNX,Xu1992} for an extensive description). \Bk

Let $u^k,\ k=0,1,2,\ldots$, be the solution sequence of the PCG algorithm. Then the following error estimate is well-known:
\[
\|u-u^k\|_A\le 2\left(\frac{\sqrt{\kappa(BA)}-1}{\sqrt{\kappa(BA)}+1}\right)^k
\|u-u^0\|_A,
\]
which implies that the PCG method converges faster with a smaller condition number $\kappa(BA)$.


\subsection{The method of parallel subspace corrections}
The method of parallel subspace corrections (PSC) provides a particular construction of the iteration \Bd operator \Bk $B$. The starting point is a suitable decomposition of $\mathcal{V}$:
\[
\mathcal{V} = \sum_{i=0}^J \mathcal{V}_i,
\]
where $\mathcal{V}_i$ are subspaces of $\mathcal{V}$. The model problem \eqref{model_prob} can be split into sub-problems in
each $\mathcal{V}_i$ with smaller size. Throughout this paper, we use the following operators, for $i=0,1,\ldots,J$:
\begin{itemize}
\item $Q_i : \mathcal{V} \rightarrow \mathcal{V}_i$ the projection in \Bd the inner product $(\cdot,\cdot):=(\cdot,\cdot)_{L^2(\Omega)}$; \Bk
\item $I_i : \mathcal{V}_i \rightarrow \mathcal{V}$ the natural inclusion;
\item \Bd $P_i : \mathcal{V} \rightarrow \mathcal{V}_i$ the projection in the inner product $(\cdot,\cdot)_A:=(A\cdot,\cdot)$; \Bk
\item $A_i : \mathcal{V}_i \rightarrow \mathcal{V}_i$ the restriction of $A$ to the subspace $\mathcal{V}_i$;
\item $R_i : \mathcal{V}_i \rightarrow \mathcal{V}_i$ an approximation of $A_i^{-1}$.
\end{itemize}

This method performs the correction on each subspace in parallel, {\Bd with the operator $B$ defined} by
\begin{equation}\label{def_B}
B:=\sum_{i=0}^{J} I_i R_i I^t_i \quad \big(= \sum_{i=0}^{J} R_i Q_i\big),
\end{equation}
\Bd where $I_i^t$ denote the adjoint of $I_i$ with respect to $(\cdot,\cdot)$. It is readily to check that $Q_iA=A_iP_i$ and $Q_i=I_i^t$ with $(I^t_i u ,v_i)=(u,I_iv_i)$. \Bk 

The convergence analysis of (PSC) is based on the following two important properties \cite{CNX}:

\noindent {\bf (A1) Stable Decomposition} For any $v\in\mathcal{V}$, there exists a decomposition $v=\sum_{i=0}^J
v_i,\ v_i\in\mathcal{V}_i,\ i=0,\ldots,J$ such that
\begin{equation*}
\sum_{i=0}^{J} \|v_i\|_A^2 \le C_1 \|v\|_A^2.
\end{equation*}

\noindent {\bf (A2) Strengthened Cauchy-Schwarz (SCS) inequality} For any $u_i,v_i\in\mathcal{V}_i,
i=0,\ldots,J$
\begin{equation*}
\left| \sum_{i=0}^{J} \sum_{j=i+1}^{J} (u_i,v_j)_A \right| \le C_2
\left(\sum_{i=0}^{J} \|u_i\|_A^2 \right)^{1/2}
\left(\sum_{i=0}^{J} \|v_i\|_A^2 \right)^{1/2}.
\end{equation*}

For a space decomposition satisfying both properties, one can prove the following result on the preconditioned linear system \cite{CNX}:
\begin{theo} \label{BPX_preconditioning}
Let $\mathcal{V}=\sum_{i=0}^J \mathcal{V}_i$ be a space decomposition satisfying \Aone and \Atwo,
and let $R_i$ be SPDs for $i=0,\ldots,J$ such that
\begin{equation}\label{assumption_smoother}
C_4^{-1}\|u_i\|_A^2\le (R_i^{-1}u_i,u_i) \le C_3 \|u_i\|_A^2 {\Bd \quad \text{ for any } u_i \in {\cal V}_i.}
\end{equation}
Then $B$ defined by \eqref{def_B} is SPD and
\begin{equation} \label{eq:cond}
\kappa(BA) \le (1+2C_2) C_1 C_3 C_4.
\end{equation}
\end{theo}

The goal of the paper is the construction of a preconditioner like \eqref{def_B} for T-splines, and the proof that \Aone and \Atwo\ are satisfied.

\section{Splines}
In this section we recall the definition and main properties of B-splines, mainly to fix the notation. For a more extensive description on the use of splines in isogeometric analysis, the reader is referred to \cite{Hughes_Cottrell_Bazilevs,IGA-book}, and to \cite{IGA-acta} for a mathematical perspective.

\subsection{Univariate B-splines}

Given two positive integers $p$ and $n$, we say that $\Xi:=\{\xi_0,\ldots,\xi_{n+p}\}$ is a $p$-open knot vector if
\begin{equation*}
0=\xi_0=\cdots=\xi_{p}<\xi_{p+1}\le \cdots \le \xi_{n-1} < \xi_{n}=\cdots=\xi_{n+p}=1,
\end{equation*}
where repeated knots are allowed. From the knot vector $\Xi$, univariate B-spline basis functions of degree $p$ are defined recursively using the Cox-De Boor formula (see \cite{DeBoor}). The definition of each B-spline ${B}_{i,p}, \ i = 0, \ldots, n-1$, is determined by a $p+2$ local knot vector $\Xi_{i,p}=\{\xi_i,\ldots,\xi_{i+p+1}\}$. Whenever necessary, we will stress it by adopting the following notation:
\begin{equation*}
{B}_{i,p}(\zeta)={B}[\Xi_{i,p}](\zeta), \quad \zeta \in (0,1).
\end{equation*}
Thus, the basis function ${B}_{i,p}$ has support
\[
\mbox{supp}({B}_{i,p})=[\xi_i,\xi_{i+p+1}].
\]
Throughout the paper, we assume that the maximum multiplicity of the internal knots is less than or equal to the degree $p$, that is, the B-spline functions are at least continuous.

Let us select from $\Xi$ a subset $\{\xi_{i_k}, k=0,\ldots,N\}$
of non-repeated knots, or breakpoints, with $\xi_{i_0}=0$, $\xi_{i_N}=1$. We point out that the local mesh
size of the element $I_k=(\xi_{i_{k}},\xi_{i_{k+1}})$ is called $h_k=\xi_{i_{k+1}}-\xi_{i_k} $ for $k=0,\ldots,N-1$.
Moreover, to the element $I_k=(\xi_{i_k},\xi_{i_{k+1}})$, that can be written as
$(\xi_{j},\xi_{j+1})$ for a certain $j$, we associate the {\it support extension}
$\widetilde{I}_j$ defined by
\begin{equation}\label{suppext_univariate}
\widetilde{I}_j:=(\xi_{j-p},\xi_{j+p+1}).
\end{equation}

We denote by
\begin{equation*}
\mathcal{S}_p(\Xi):= \mbox{span}\{{B}_{i,p}, i=0,\ldots,n-1\}.
\end{equation*}
the univariate B-splines space spanned by those B-splines of degree $p$. The functions ${B}_{i,p}$ are a \emph{partition of unity}, as shown in \cite{Schumi}.





Following  \cite[Theorem~4.41]{Schumi},  we define suitable
functionals $\lambda_{j,p} =
\lambda[\xi_{j},\ldots,\xi_{j+p+1}]$, for $0
\le j \le n-1 $, which are dual to the B-splines basis functions, that is
%
\begin{equation*}
\lambda_{j,p}({B}_{i,p})=\delta_{ij},\qquad 0\le i,j\le n-1,
\end{equation*}
where $\delta_{ij}$ is the Kronecker delta. The following estimate of the functionals
$\lambda_{j,p}$ will be useful in the sequel.
\begin{lemm}\label{lemm1}
If $f\in L^q(\xi_j,\xi_{j+p+1})$, with $1 \le q \le +\infty $
then \[ |\lambda_{j,p}(f)|\lec |\xi_{j+p+1} -
\xi_j|^{-1/q}\|f\|_{L^q(\xi_j,\xi_{j+p+1})}.
\]
\end{lemm}
{\it Proof.} The proof can be found in \cite[Theorem~4.41]{Schumi}. \hfill $\square$

We note that these dual functionals are locally defined and depend only on the corresponding local knot vector, namely,
\begin{equation*}
\lambda_{i,p}(f) = \lambda[\Xi_{i,p}](f).
\end{equation*}

Let $\Pi_{p,\Xi}$ be the projection
that is built from the dual basis as detailed in \cite[Theorem~12.6]{Schumi}, that is,
\begin{equation}\label{proj_univariate}
\Pi_{p,\Xi}: L^2([0,1]) \rightarrow \mathcal{S}_p(\Xi),\qquad \Pi_{p,\Xi}(f)=\sum_{j=0}^{n-1}
\lambda_{j,p}(f){B}_{j,p}.
\end{equation}

\begin{assumpt}\label{assumpt_quasiuniform}
$\{\xi_{i_0},\xi_{i_1},\ldots,\xi_{i_N}\}$ is locally
quasi-uniform, that is, there exists a constant $\theta \ge 1$ such that the
mesh sizes $h_k = \xi_{i_{k+1}}-\xi_{i_k}$ satisfy the relation $\theta^{-1} \le h_k /h_{k+1} \le \theta$, for $ k = 0, \ldots , N-2$.
\end{assumpt}

\begin{prop} For any non-empty knot span $I_k=(\xi_{i_k},\xi_{i_{k+1}})$, we have
\begin{equation*}
\|\Pi_{p,\Xi}(f)\|_{L^2(I_k)} \le C \| f \|_{L^2(\widetilde{I}_k)},
\end{equation*}
where the constant $C$ depends only on the degree $p$. Moreover, if Assumption \ref{assumpt_quasiuniform} holds, we have
\begin{equation*}
|\Pi_{p,\Xi}(f)|_{H^1(I_k)} \le C | f |_{H^1(\widetilde{I}_k)},
\end{equation*}
where the constant $C$ depends only on $p$ and $\theta$ and $H^1$ denotes the classical Sobolev norm.
\end{prop}
{\it Proof.} The proof can be found in \cite[Proposition~2.2]{IGA-acta}. \hfill
$\square$

\subsection{Multivariate splines}

Multivariate B-splines can be constructed by means of tensor products. We discuss here the bivariate case, the higher-dimensional case being analogous.

For $d=1,2$, assume that $n_d \in \mathbb{N}$, the degree $p_d$ and the $p_d$-open knot vector
$\Xi_d=\{\xi_{d,0},\ldots,\xi_{d,n_d+p_d}\}$ are given. We set the polynomial degree vector
${\mathbf p}=(p_1,p_2)$ and ${\bf \Xi}=\{\Xi_1, \Xi_2\}$. We introduce a set of multi-indices ${\bf I}=\{{\bf i}=(i_1,i_2):0\le i_d\le n_d-1\}$ and for each multi-index ${\bf i}=(i_1,i_2)$, we define the local knot vector
\[
{\bf \Xi}_{\bf i, p} = \{ \Xi_{i_1,p_1}, \Xi_{i_2,p_2} \}.
\]
Then we can introduce the set of multivariate B-splines
\begin{equation*}
\left\{ {B}_{\bf i,p}({\bzeta}) = {B}[\Xi_{i_1,p_1}](\zeta_1){B}[\Xi_{i_2,p_2}](\zeta_2) , \quad \mbox{for\ all\ } {\bf i} \in {\bf I}\right\}.
\end{equation*}
The spline space in the parametric domain ${\Omega}=(0,1)^2$ is then
\[
S_{\bf p}({\bf \Xi})=\mbox{span}\{{B}_{\bf i,p}({\bzeta}),\quad {\bf i} \in {\bf I}\}.
\]

We also introduce the set of non-repeated interface knots for each direction $\{\xi_{d, i_0},\ldots,\xi_{d,i_{N_d}}\}$, $d = 1,2$, which determine the intervals $I_{d,j_d}=(\xi_{d,i_{j_d}}, \xi_{d,i_{j_d+1}})$, for $0 \le j_d \le N_{d}-1$. These intervals lead to the Cartesian grid ${\mathcal{M^B}}$ (or simply ${\cal M}$) in the unit domain ${\Omega} = (0,1)^2$, also called the B\'ezier mesh:
\begin{equation*}
\mathcal{M^B}=\{Q_{\bf j}=I_{1,j_1} \times I_{2,j_2}, \ \mbox{for\ } 0 \le j_d \le N_d-1 \}.
\end{equation*}
For a generic element $Q_{\bf j}$, we also define its support extension
\[
\widetilde{Q}_{\bf j} = \widetilde{I}_{1,j_1}\times \widetilde{I}_{2,j_2},
\]
with $\widetilde{I}_{d,j_d}$ the univariate support extension by \eqref{suppext_univariate}.


\subsubsection{Multivariate quasi-interpolants}
We note that when the univariate quasi-interpolants are defined from a dual basis, as in \eqref{proj_univariate}, then the multivariate quasi-interpolant is also defined from a dual basis. Indeed,
we have
\begin{equation}\label{multi_quasiint}
{\bf \Pi}_{\bf p,\Xi}(f) = \sum_{{\bf i}\in{\bf I}} \lambda_{\bf i,p}(f){B}_{\bf i, p},
\end{equation}
where each dual functional is defined from the univariate dual bases by the expression
\begin{equation*}
\lambda_{\bf i,p} = \lambda_{i_1,p_1} \otimes \lambda_{i_2,p_2}.
\end{equation*}

\section{T-splines}\label{section4}
The main drawback of B-splines is their tensor-product structure, which prevents local refinement as required by adaptive methods. One of the alternatives is the use of T-splines \cite{Sederberg_Zheng}, a superset of B-splines that allows for local refinement. The use of T-splines in IGA was first explored in \cite{Bazilervs_Calo_Cottrell_Evans,Dorfel_Juttler_Simeon}, and it has led to a growing interest for the analysis of their mathematical properties. In this section we are collecting mathematical results from \cite{LZSHS12,BBCS12,BBSV2013} (linear-independence, dual basis and projectors), \cite{Li-Scott,Bressan_Buffa_Sangalli} (nestedness and space characterization) and \cite[Section 7]{IGA-acta} (local linear independence), following mainly the notation in \cite{IGA-acta}.

We restrict ourselves to T-splines where refinement is always performed by bisection, and the multiplicity is never reduced. A more general setting can be considered, but it would only add technical difficulties without adding more insight.

\subsection{T-mesh defined by bisection} \label{sec:T-bisection}
As in the previous section, let us assume that we are given the degrees $p_d$, the integers $n_d$ and the open knot vectors $\Xi_d$, and let us denote $m_d = n_d + p_d$ for $d = 1,2$. For simplicity we assume that the internal knots in $\Xi_d$ are not repeated and \emph{equally spaced}, so the element size in each parametric direction can be denoted by $h_1$ and $h_2$. Our starting point is the index mesh ${\cal T}_0$ in the index domain $[0,m_1]\times [0,m_2]$, defined as the Cartesian grid of unit squares
$$
{\cal T}_0 = \{[j_1, j_1+1] \times [j_2, j_2+1] : j_1 = 0, \ldots, m_1-1; \, j_2= 0, \ldots, m_2-1\},
$$
which is associated to the tensor-product B-spline space $S_\bp(\boldsymbol \Xi)$.

To define T-splines by bisection, we start introducing, for any integer $\ell \ge 0 $ and for $d = 1,2$, the set of rational indices
\footnotesize
\begin{equation*}
{\cal I}^\ell_d = \left\{0, \ldots, p_d, p_d + \frac{1}{2^\ell},\ldots, p_d + \frac{2^\ell-1}{2^\ell}, p_{d}+1,\ldots, n_{d}-1, n_{d}-1 + \frac{1}{2^\ell},\ldots, n_{d}-1 + \frac{2^\ell-1}{2^\ell}, n_d, n_d+1, \ldots, n_d + p_d \right\},
\end{equation*}
\normalsize
and we notice that ${\cal I}^\ell_d \subset {\cal I}^{\ell'}_d$ if $\ell \le \ell'$. We also define the ordered knot vectors
\[
\Xi^\ell_d = \{\xi_{d,k}, \, k \in {\cal I}^\ell_d \} = \left\{ \xi_{d,0}, \xi_{d,1}, \ldots, \xi_{d,p_d}, \xi_{d,p_d+\frac{1}{2^\ell}}, \ldots, \xi_{d,n_{d}-1 + \frac{2^\ell-1}{2^\ell}}, \xi_{d,n_d}, \xi_{d,n_d+1}, \ldots, \xi_{d,n_d+p_d} \right\},
\]
in a recursive way: starting from $\Xi_d^0 = \Xi_d$, for $\ell > 0$ and for any new index $k \in {\cal I}^\ell_d \setminus {\cal I}^{\ell-1}_d$, we define the knot
\[
\xi_{d,k} = \frac{1}{2} (\xi_{d,k-\frac{1}{2^\ell}} + \xi_{d,k+\frac{1}{2^\ell}}),
\]
which is well defined because $k - \frac{1}{2^\ell}, k + \frac{1}{2^\ell} \in {\cal I}_d^{\ell-1}$. Clearly, $\Xi_d^\ell \subset \Xi_d^{\ell+1}$ for $\ell \ge 0$, and the interval size is $h_{d,\ell} = h_d/2^\ell$. Notice that in this procedure we do not introduce new knots between the repeated knots of the open knot vector.

We also define, for an arbitrary rectangular element in the index domain
$\tau = [k_1, k_1 + t_1] \times [k_2, k_2 + t_2]$, with indices $k_d, k_d + t_d \in {\cal I}^\ell_d$, the following bisection operators (see \cite[Definition~2.5]{Morgenstern_Peterseim}):
\begin{equation*}
\begin{array}{l}
\bisect_x(\tau) =
\left \{
\begin{array}{ll}
\{ [k_1, k_1 + t_1/2] \times [k_2, k_2 + t_2], [k_1 + t_1/2, k_1 + t_1] \times, [k_2, k_2 + t_2] \} & \text{ if } \xi_{1,k_1} \ne \xi_{1,k_1+t_1}, \\
\tau & \text{ if } \xi_{1,k_1} = \xi_{1,k_1+t_1},
\end{array}
\right. \\
\bisect_y(\tau) =
\left \{
\begin{array}{ll}
\{ [k_1, k_1 + t_1] \times [k_2, k_2 + t_2/2], [k_1, k_1 + t_1] \times, [k_2 + t_2/2, k_2 + t_2] \} & \text{ if } \xi_{2,k_2} \ne \xi_{2,k_2+t_2},\\
\tau & \text{ if } \xi_{2,k_2} = \xi_{2,k_2+t_2}.
\end{array}
\right. \\
\end{array}
\end{equation*}
Notice that the bisection operators split the element in two adding a vertical and a horizontal edge, respectively, only when the corresponding length in the parametric domain is greater than zero, that is, when $\xi_{d,k_d + t_d} - \xi_{d,k_d} > 0$.

Starting from the Cartesian grid ${\cal T}_0$, we define a T-mesh ${\cal T} = {\cal T}_N$ by successive applying bisection of elements, in the form
\begin{equation} \label{eq:T-by-bisection}
{\cal T}_{k+1} = {\cal T}_k + b_{\tau_k}, \quad \tau_k \in {\cal T}_k, \; k = 0, \ldots, N-1,
\end{equation}
where we use the formal addition (see \cite[Definition~2.6]{Morgenstern_Peterseim} and also \cite[Section~3]{CNX})
\begin{equation} \label{eq:formal_addition}
{\cal T}_k + b_{\tau_k} = {\cal T}_k \setminus \{\tau_k\} \cup \bisect(\tau_k),
\end{equation}
and the bisection operator can be either $\bisect_x$ or $\bisect_y$. Moreover, we define the finest level $\Lambda$ as the minimum integer such that the $d$th coordinate of all vertices in the T-mesh ${\cal T}$ belongs to ${\cal I}_d^\Lambda$, for $d=1,2$.

\begin{remark}
We notice that a T-mesh defined with this procedure automatically satisfies the first condition in the definition of \emph{admissible mesh} in \cite[Definition~7.10]{IGA-acta}, that is, the first $p_d+1$ lines closer to each boundary are completely contained in the mesh. Instead, the second condition of not having T-junctions in the so-called frame region is not satisfied, because T-junctions may appear on the interface between the frame and the active region. In any case, these T-junctions do not affect the definition of the T-spline functions, and the mesh can still be considered ``admissible''.
\end{remark}

\subsection{Analysis suitable T-splines}
\begin{figure}[ht!]
\includegraphics[width=0.49\textwidth, trim = 3cm 0cm 3cm 0cm, clip]{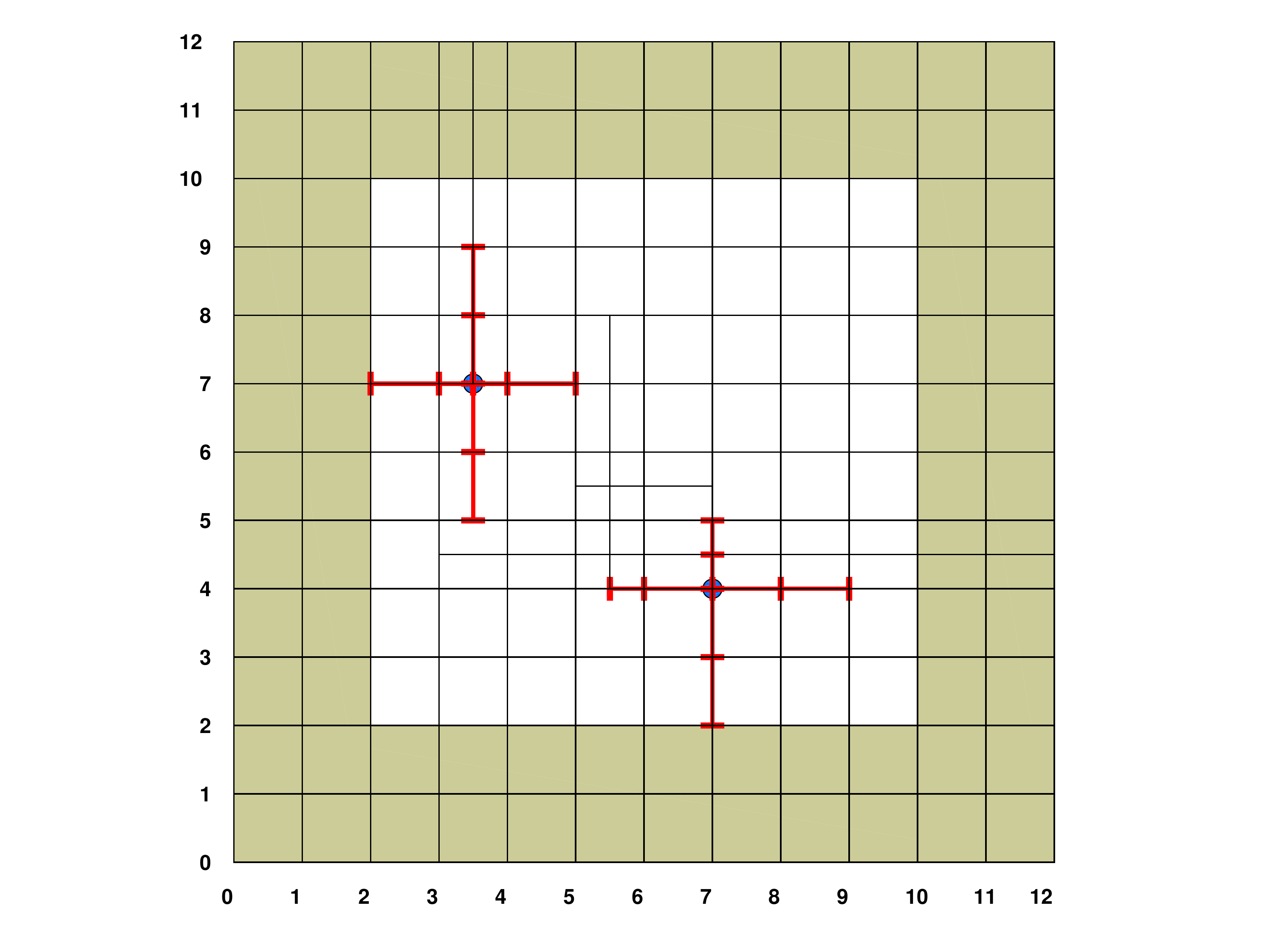}
\includegraphics[width=0.49\textwidth, trim = 3cm 0cm 3cm 0cm, clip]{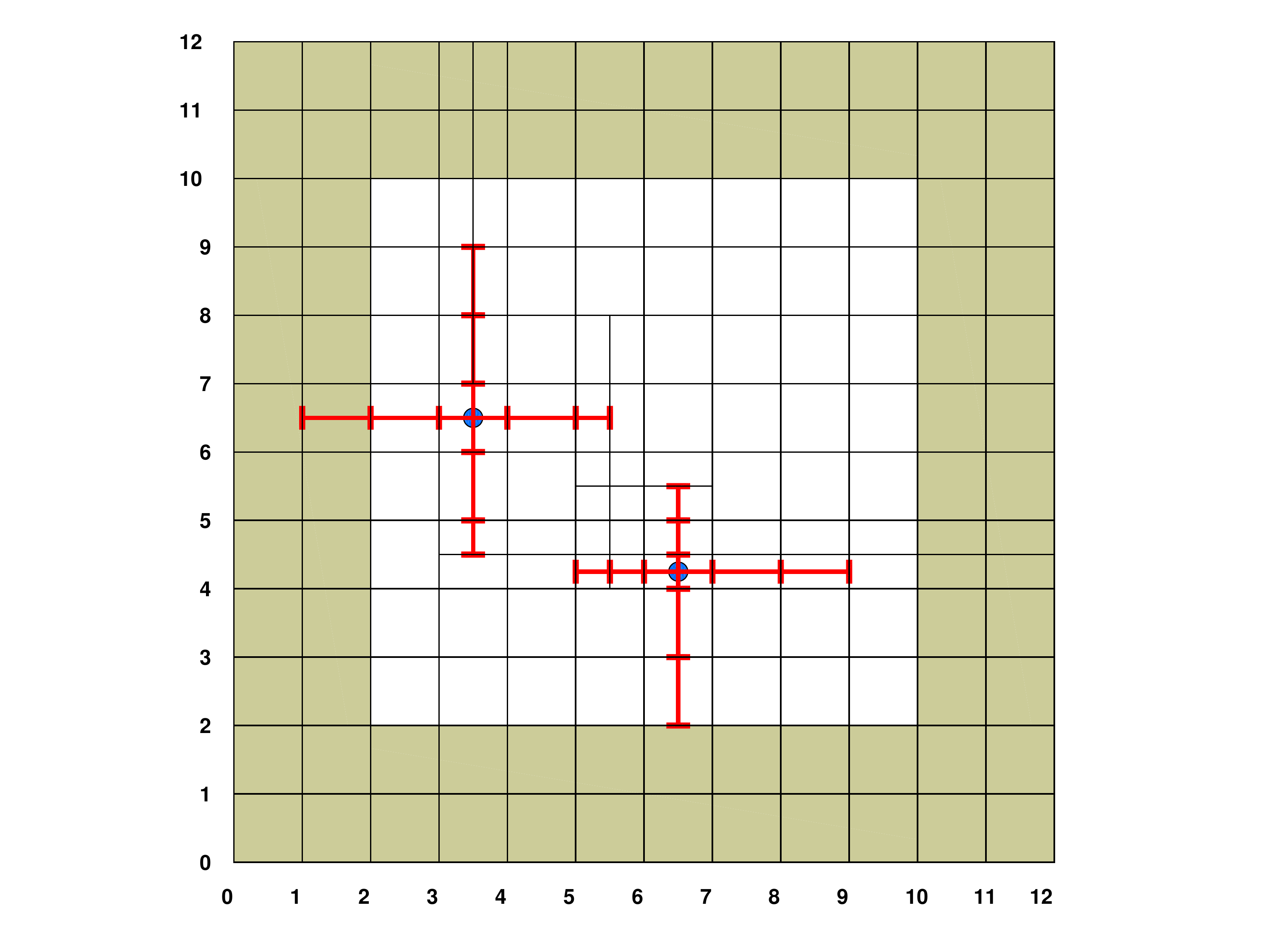}
\caption{Computation of the index vector for two basis functions, for bicubic (left) and biquartic (right) T-splines} \label{fig:anchors}
\end{figure}

To construct the blending functions associated to a T-mesh ${\cal T}$, we have
to define first the set of \emph{anchors}, that we denote by ${\cal A}_{\bf
p}({\cal T})$, and that depends on the degree, see \cite[Definition~7.13]{IGA-acta}.
These are either the set of vertices ($p_1$ and $p_2$ odd), elements ($p_1$ and
$p_2$ even), horizontal edges ($p_1$ even, $p_2$ odd) or vertical edges ($p_1$
odd, $p_2$ even) in the active region, which is the rectangle
\[
[ \lceil p_1/2 \rceil,  m_1 - \lceil p_1/2 \rceil] \times [ \lceil p_2/2 \rceil,  m_2 - \lceil p_2/2 \rceil],
\]
see the examples in Figure~\ref{fig:anchors}. We also define the horizontal (resp. vertical) skeleton of the mesh, and denote it by $\hsk({\cal T})$ (resp. $\vsk({\cal T})$), as the union of all horizontal (resp. vertical) edges and vertices. The union of $\hsk({\cal T})$ and $\vsk({\cal T})$ will be called skeleton. Then, for each anchor we construct an ordered horizontal index vector of $p_1+2$ indices, tracing a horizontal line from the anchor and collecting the closest $\lfloor (p_1+2)/2 \rfloor$ indices leftwards and rightwards where the line intersects the vertical skeleton of the mesh, plus the index of the anchor if the degree is odd, see \cite[Definition~7.14]{IGA-acta} and Figure~\ref{fig:anchors}. A vertical index vector of $p_2+2$ indices is constructed in an analogous way tracing a vertical line passing through the anchor.

Then, given an anchor $\bA \in {\cal A}_\bp({\cal T})$ with horizontal and vertical index vectors
\begin{equation} \label{eq:index-vectors}
hv_{\bp}(\bA) = \{i_1, \ldots, i_{p_1+2} \} \subset {\cal I}^\Lambda_1, \quad
vv_{\bp}(\bA) = \{j_1, \ldots, j_{p_2+2} \} \subset {\cal I}^\Lambda_2,
\end{equation}
we define the local knot vectors
\[
\Xi_{\bA,1,p_1} = \{\xi_{1,i_1}, \ldots, \xi_{1,i_{p_1+2}} \} \subset \Xi_1^\Lambda,
\quad \text{ and } \; \Xi_{\bA,2,p_2} = \{\xi_{2,j_1}, \ldots, \xi_{2,j_{p_2+2}}
\} \subset \Xi_2^\Lambda,
\]
and from these local knot vectors we define the associated bivariate function
\begin{equation*}
B_{\bA,\bp}(\bzeta) =
B[\Xi_{\bA,1,p_1}](\zeta_1) B[\Xi_{\bA,2,p_2}] (\zeta_2).
\end{equation*}
We will denote by $S({\cal A}_\bp({\cal T})) = {\rm span} \{B_{\bA,\bp} : \bA \in {\cal A}_\bp({\cal T}) \}$ the space generated by the T-splines.

A sufficient condition to guarantee linear independence of the T-spline functions is the dual compatibility condition, introduced in \cite{BBCS12,BBSV2013} and equivalent to the analysis-suitability condition in \cite{Li-Scott}. The former was generalized to arbitrary dimension in \cite{IGA-acta}. We say that two local knot vectors \emph{overlap} if both of them can be written as subvectors, with consecutive indices, of the same global knot vector. Then we say that the T-mesh is \emph{dual compatible} if for each pair of anchors $\bA', \bA'' \in {\cal A}_\bp({\cal T})$, with $\bA' \ne \bA''$, there exists a direction $d \in \{1,2\}$ such that the local knot vectors $\Xi_{\bA',d,p_d}$ and $\Xi_{\bA'',d,p_d}$ are different and overlap. See \cite[Section~7.1]{IGA-acta} for details.

Assuming that the T-mesh is dual compatible, then the functionals
\begin{equation*}
\{\lambda_{\bA,\bp}, \bA \in {\cal A}_\bp({\cal T}) \}, \quad \lambda_{\bA,\bp} := \lambda[\Xi_{\bA,1,p_1}] \otimes \lambda[\Xi_{\bA,2,p_2}],
\end{equation*}
form a dual basis for $\{{B}_{\bf A,p}: {\bf A}\in \mathcal{A}_{\bf p}(\mathcal{T}) \}$ (see \cite[Proposition~7.3]{IGA-acta}). Moreover, we can build the projection operator ${\bf \Pi}_\bp^{\cal T} \equiv {\bf \Pi}_{\mathcal{A}_\bp(\mathcal{T})}:
L^2({\Omega}) \rightarrow S(\mathcal{A}_{\bf p}(\mathcal{T}))$ by
\begin{equation}\label{Projection_AST}
{\bf \Pi^{\cal T}_{\bf p}}(f)({\bzeta})
 := \sum_{\bA \in {\cal A}_\bp({\cal T})} \lambda_{{\bf A,p}}(f){B}_{\bf A, p}({\bzeta})\quad
\mbox{\ for\ all\ } f\in L^2({\Omega}),\ \mbox{\ and\ } {\bzeta} \in {\Omega}.
\end{equation}

\begin{figure}[ht!]
\centering
\begin{subfigure}[Extended T-mesh]{
\includegraphics[width=0.47\textwidth,trim=2cm 1cm 2cm 0cm, clip]{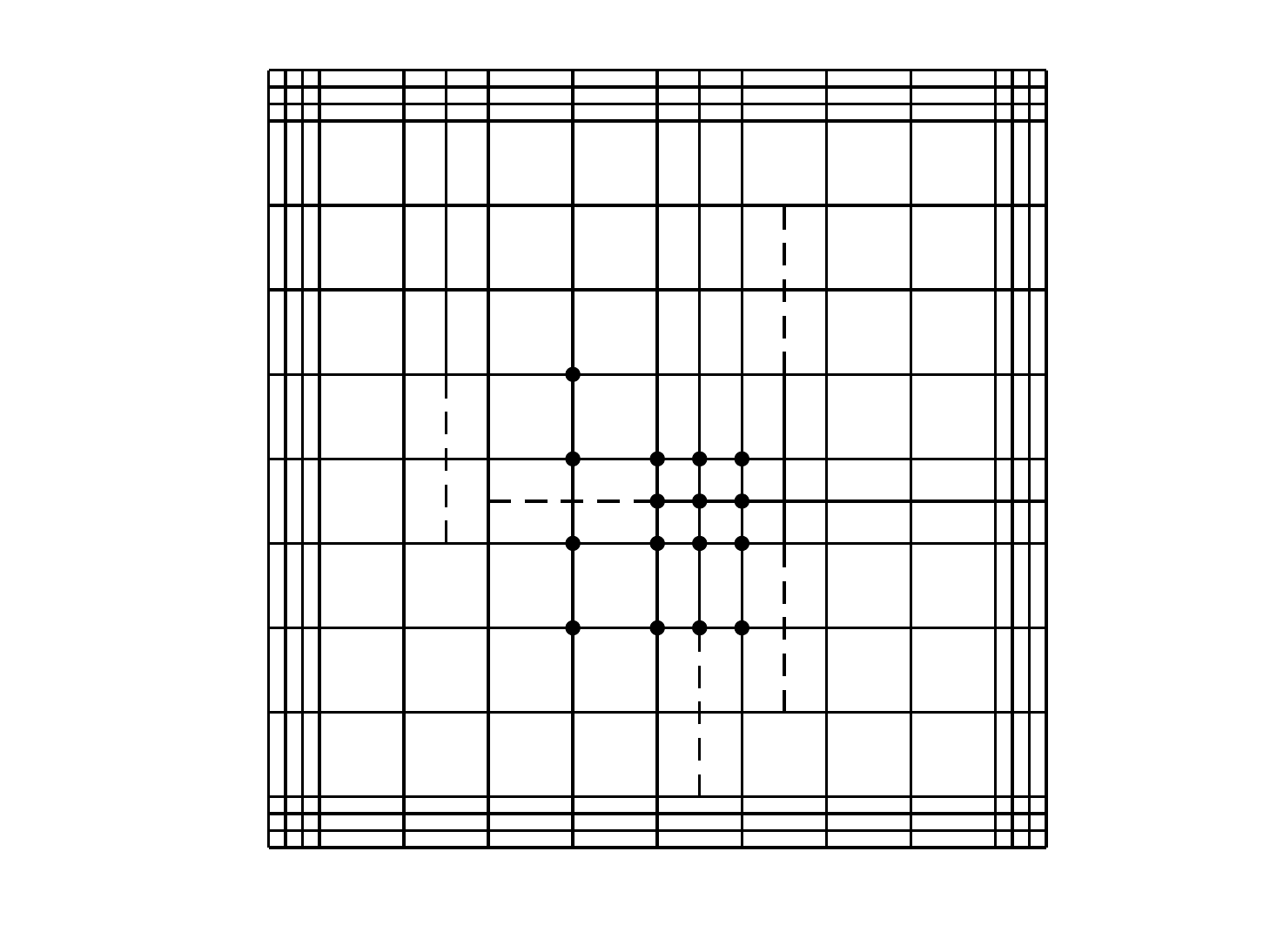}
}
\end{subfigure}
\begin{subfigure}[The corresponding B\'ezier mesh]{
\includegraphics[width=0.47\textwidth,trim=2cm 1cm 2cm 0cm, clip]{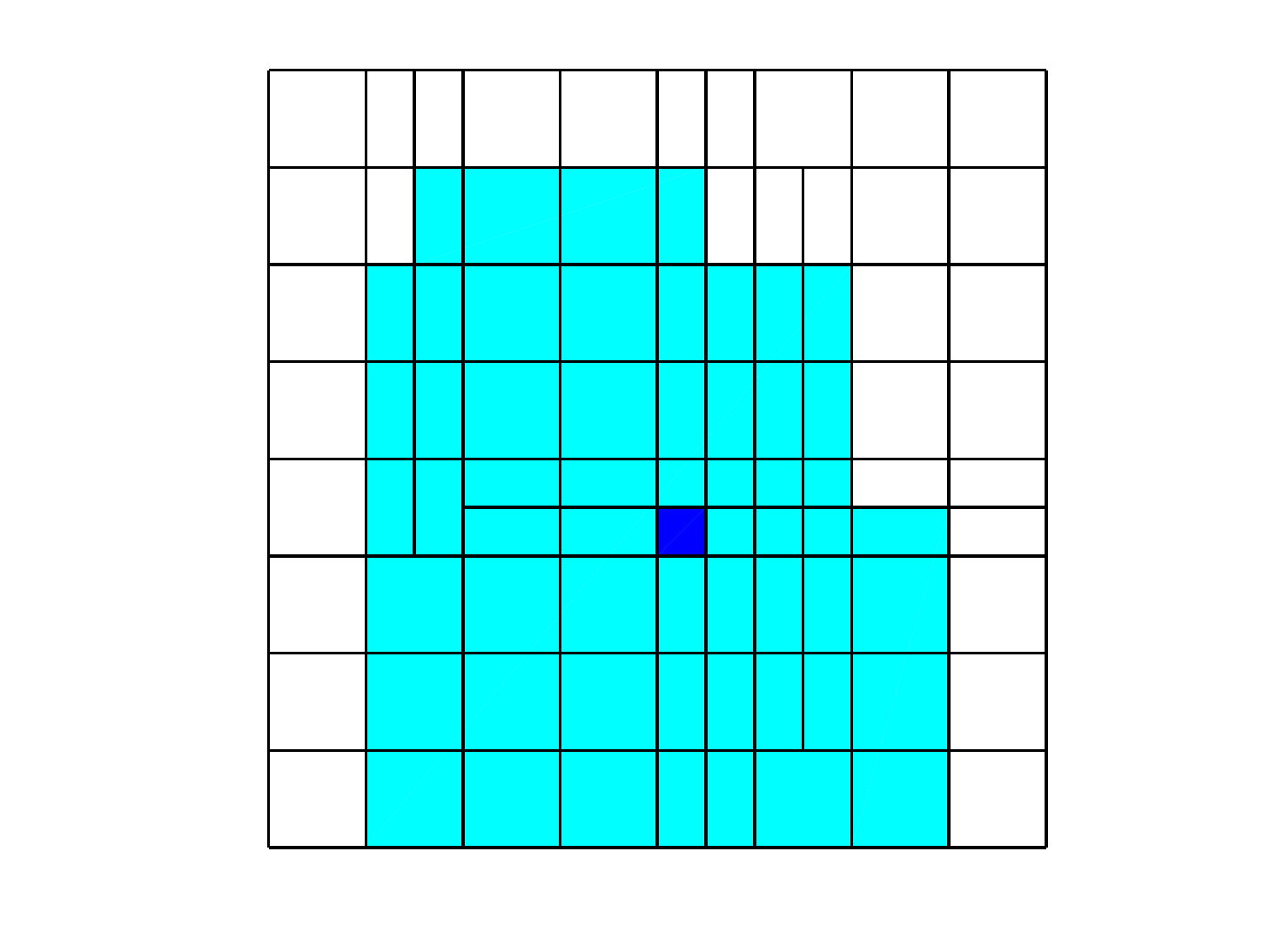}
}
\end{subfigure}
\caption{An extended T-mesh and the corresponding B\'ezier mesh for degree ${\bf p}=(3,3)$. Repeated knots in the T-mesh, which lead to empty elements, are represented with small separation between lines. In the B\'ezier mesh, we highlight a generic element $Q$ (in blue) and its support extension $\widetilde{Q}$ (in sky blue). The functions that do not vanish in $Q$ are marked on the left by its anchors} \label{fig:beziermesh_extsupport}
\end{figure}

For the analysis of the projector properties, we make use of the B\'ezier mesh,
that we define as in \cite[Section~7.3]{IGA-acta}. The B\'ezier mesh is different
from the T-mesh, and plays a similar role to the mesh in finite elements. We
start recalling that T-junctions are internal vertices of the T-mesh with
valence equal to three, that can be grouped in horizontal ($\vdash,\dashv$) and
vertical ($\top, \bot$) T-junctions. For a T-junction of type $\dashv$ with
index coordinates $(\bar \imath, \bar \jmath)$, we define the extension as the
minimal horizontal line that, passing through the T-junction, intersects
$\lfloor p_1/2 \rfloor$ (closed) vertical edges to its left, and $\lceil p_1/2
\rceil$ to its right. The extensions in the other cases are defined in a similar
way, using $p_2$ for vertical T-junctions. Then, we define ${\rm ext}_\bp ({\cal
T})$ the \emph{extended T-mesh} in the index domain adding to ${\cal T}$ all the
T-junction extensions. The B\'ezier mesh ${\cal M^T}$ (or simply ${\cal M}$)
associated to ${\cal T}$ is then defined as the collection of \emph{non-empty
elements} in the domain $\Omega$ of the form
\begin{equation*}
Q = (\xi_{1,i_1}, \xi_{1,i_2}) \times (\xi_{2, j_1}, \xi_{2, j_2}) \ne \emptyset, \quad \text{ with } (i_1, i_2) \times (j_1,j_2) \in {\rm ext}_\bp ({\cal T}),
\end{equation*}
as depicted in Figure~\ref{fig:beziermesh_extsupport}.

For a B\'ezier element $Q \in {\cal M^T}$, and in general for any subdomain $Q \subset \Omega$, we define the support extension as
the union of the supports of the functions whose support intersects $Q$, that
is
\begin{equation}\label{supportext_Q}
\tilde Q := \bigcup_{\bA \in {\cal A}_Q} {\rm supp} (B_{\bA,\bp}), \text{ with } {\cal A}_Q = \{ \bA \in {\cal A}_\bp({\cal T}) : {\rm int}({\rm supp} (B_{\bA,\bp})) \cap Q \ne \emptyset \}.
\end{equation}
where int($C$) denotes the interior of a set $C$.
Moreover, we define $\bar Q$ as the smallest rectangle in $\Omega$ containing $\tilde Q$. The following result holds (see \cite[Proposition~7.7]{IGA-acta}):
\begin{prop}\label{L2_stable_Pi_AST}
Let $\mathcal{T}$ be a dual compatible T-mesh. Then there exists a constant $C$,
depending only on ${\bf p}$, such that for any B\'ezier element $Q \in {\cal
M^T}$ the projector \eqref{Projection_AST} satisfies
\[
\|{\bf \Pi}^{\cal T}_{\bf p}(f)\|_{L^2(Q)} \le C \| f \|_{L^2(\widetilde{Q})},\ \mbox{for\ all\ } f\in L^2({\Omega}).
\]
\end{prop}

Finally, we notice that for each anchor the index vectors
\eqref{eq:index-vectors} define a local Cartesian grid of $(p_1+1)
(p_2+1)$ cells, called \emph{tiled floor} in \cite{BBCS12,BBSV2013}. Moreover,
we also define the \emph{parametric tiled floor}, as the set
of non-empty cells
\[
[\xi_{1,i_k}, \xi_{1,i_{k+1}}] \times [\xi_{2,j_{k'}}, \xi_{2,j_{k'+1}}] \ne
\emptyset, \text{ with } \,  i_k,i_{k+1} \in hv_\bp(\bA), \; j_{k'}, j_{k'+1}
\in vv_\bp(\bA).
\]
We remark that in
general, the cells of the tiled floor do not coincide with the elements of the
T-mesh, and the cells of the parametric tiled floor do not coincide
with the elements of the B\'ezier mesh.

\subsection{Analysis suitable T-splines by bisection} \label{sec:ASTS-bisection}
In order to apply BPX preconditioners to analysis-suitable T-splines, it is necessary to define a suitable refinement procedure that provides a multilevel structure. In this section we adopt the refinement strategy introduced and analysed in \cite{Morgenstern_Peterseim}, and present a new local quasi-uniformity result necessary for the analysis of multilevel preconditioners. The idea in that work is to refine by bisection alternating the refinement direction, and whenever a new edge is added, a recursive algorithm is called to refine the elements in the neighborhood, ensuring that the condition of analysis-suitability is preserved. One of the advantages of this refinement algorithm is that it allows to associate a level (or generation) to each element and function, as required by multilevel methods.

We start setting the \emph{generation} $g(\tau) = 0$ for all the elements of the Cartesian grid $\tau \in {\cal T}_0$. Then, the T-mesh ${\cal T}$ is defined as in \eqref{eq:T-by-bisection}--\eqref{eq:formal_addition}, choosing the bisection
\begin{equation*}
\bisect(\tau) = \left \{
\begin{array}{l}
\bisect_x (\tau) \quad \text{ if } g(\tau) \text{ is even},\\
\bisect_y (\tau) \quad \text{ if } g(\tau) \text{ is odd},
\end{array}
\right.
\end{equation*}
and setting the generation $g(\tau') = g(\tau) + 1$ for $\tau' \in \bisect(\tau)$ (see~\cite[Definition~2.6]{Morgenstern_Peterseim}). Moreover, we say that the bisection $b_\tau$ has generation $g(b_\tau) = g(\tau) + 1$. Without loss of generality (see also the proof of Theorem~3.6 in \cite{Morgenstern_Peterseim}), in the following we will assume that the bisections in \eqref{eq:T-by-bisection} are ordered by their generation, that is, if $k > k'$ then $g(b_{\tau_k}) \ge g(b_{\tau_{k'}})$, and the same relation holds for the elements generated by the two bisections. In the following, we will denote by $L$ the finest generation (or level), that is, the generation of the last bisection.

\begin{remark}
We recall that elements with zero length in one parametric direction are not bisected in that direction, see Section~\ref{sec:T-bisection}. However, when applying the bisection operator their generation is increased by one, in such a way that the next time the bisection operator is applied, they will be refined in the other direction.
\end{remark}

\begin{figure}[ht!]
\centering
\includegraphics[width=0.6\textwidth]{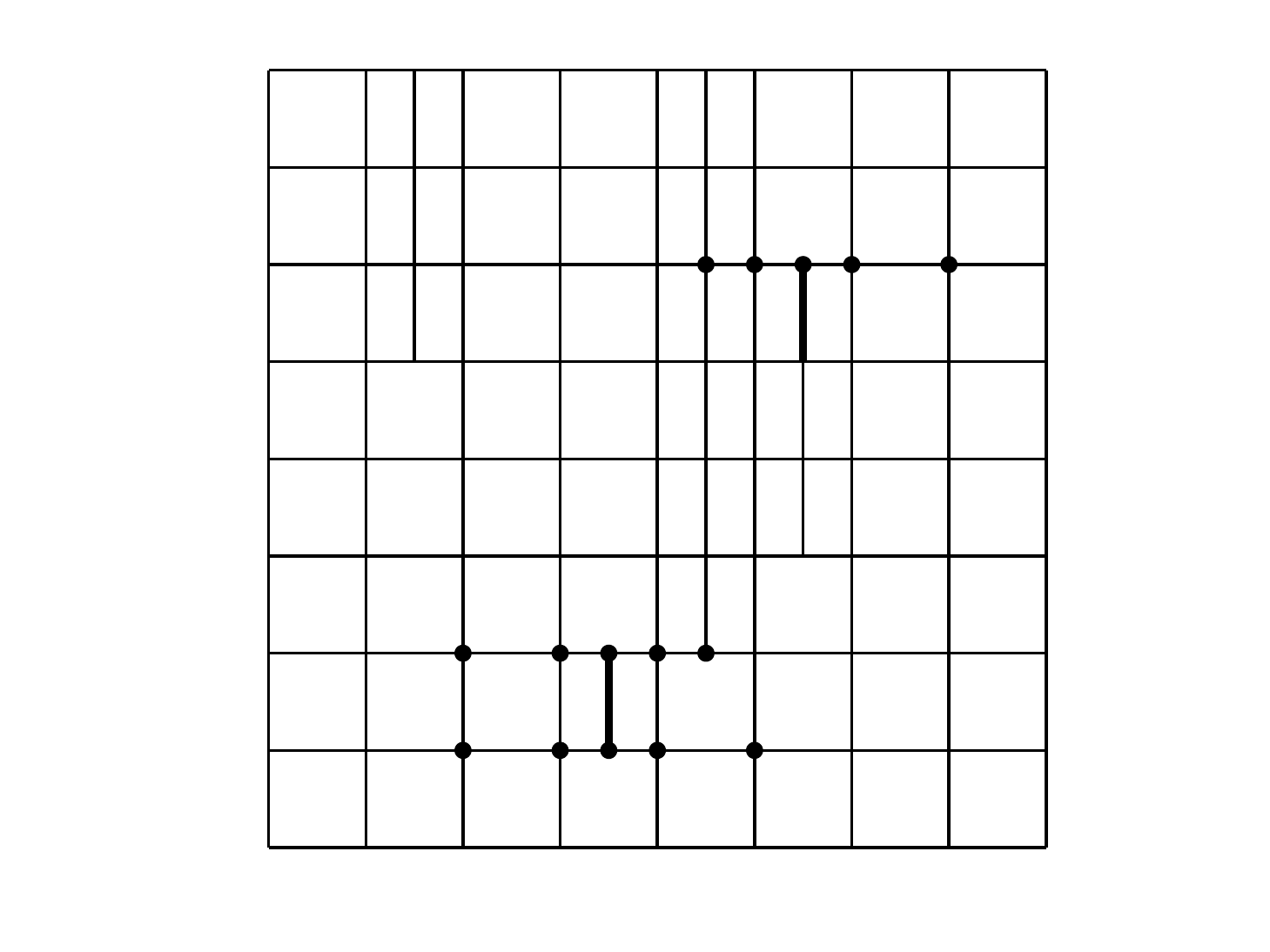}
\caption{Two bisections $b_{\tau_{k-1}}$ (lower left in bold black) and $b_{\tau_{k}}$ (upper right in bold black) of the same generation, and their respective collections $\Phi_{k-1}$ and $\Phi_{k}$ of bicubic T-spline functions where anchors near each bisection represent the associated T-splines newly appeared or modified} \label{fig:mesh_definition_Vk}
\end{figure}

To define the generation of the functions, we denote by $\Phi_0 = \{B_{\bA,\bp}: \bA \in {\cal A_\bp}({\cal T}_0) \}$ the functions of the tensor product spline space in the coarsest mesh, and then define
for $k = 1, \ldots, N$ the collection of T-spline functions newly appeared or modified after the bisection $b_{\tau_{k-1}}$ as
\[
\Phi_{k}:=\{{B}_{\bf A,p}: {\bf A}\in \mathcal{A}_{\bf p}(\mathcal{T}_{k})\} \backslash \{{B}_{\bf A,p}: {\bf A}\in \mathcal{A}_{\bf p}(\mathcal{T}_{k-1})\}.
\]
An example of the definition of $\Phi_{k}$ is shown in Figure \ref{fig:mesh_definition_Vk}. Functions in $\Phi_0$ have generation 0, whereas functions in $\Phi_{k}$ have the same generation of the bisection $b_{\tau_{k-1}}$. To alleviate notation, in the following we will denote by $\ell_k$ the generation (or level) of functions in $\Phi_k$, that is
\begin{equation*}
\ell_k := g(b_{\tau_{k-1}}) = g(\tau_{k-1}) + 1, \quad \text{ for } k = 0, \ldots, N,
\end{equation*}
with the convention that $g(\tau_{-1}) = -1$. Notice that the subscript $k$ varies from $0$ to $N$, while $\ell_k$ takes values from 0 to $L$.

To each $\Phi_{k}$, we associate a subspace
\begin{equation}\label{subspace_Vk}
\mathcal{V}_{k} :=\mbox{span}~\Phi_{k},\quad \text{ for } k = 0, \ldots, N.
\end{equation}
We also define the support of functions in $\Phi_{k}, 1\le k\le N$ and its support extension as
\begin{equation} \label{eq:omegak}
\omega_k := \bigcup_{{B}_{\bf A,p}\in\Phi_{k}} \mbox{supp}~({B}_{\bf A,p}),
\end{equation}
and
\begin{equation} \label{eq:omegak-tilde}
\tilde{\omega}_k:=  \bigcup_{\bA \in {\cal A}_{\omega_k}} {\rm supp} (B_{\bA,\bp}), \text{ with } {\cal A}_{\omega_k} := \{ \bA \in {\cal A}_\bp({\cal T}_{k}) :  {\rm int}({\rm supp} (B_{\bA,\bp})) \cap \omega_k \ne \emptyset \},
\end{equation}
respectively.

Notice that, since we alternate the directions of refinement, and recalling the notation of Section~\ref{sec:T-bisection}, the local knot vectors of a function of generation $\ell$ are contained in
\[
{\bf \Xi}^\ell = \{\Xi_1^{\lceil \ell/2 \rceil}, \Xi_2^{\lfloor \ell/2 \rfloor}\}.
\]
For convenience, we also introduce the corresponding set of rational indices, ${\cal \bf I}^\ell = \{ {\cal I}_1^{\lceil \ell/2 \rceil}, {\cal I}_2^{\lfloor \ell/2 \rfloor} \}$, and the B\'ezier mesh in the parametric domain ${\cal M}_\ell^{\cal B}$. For each generation $\ell$ we have that the mesh size is $h_\ell \approx 2^{-\ell/2}$. This important relation between generation and mesh size can be also represented using additional notation as
\[
h_\ell \approx \gamma^\ell \qquad \mbox{\ with\ } \gamma=2^{-1/2} \in (0,1).
\]

\begin{figure}[ht!]
\includegraphics[width=0.5\textwidth,trim=1cm 1cm 1cm 0cm, clip]{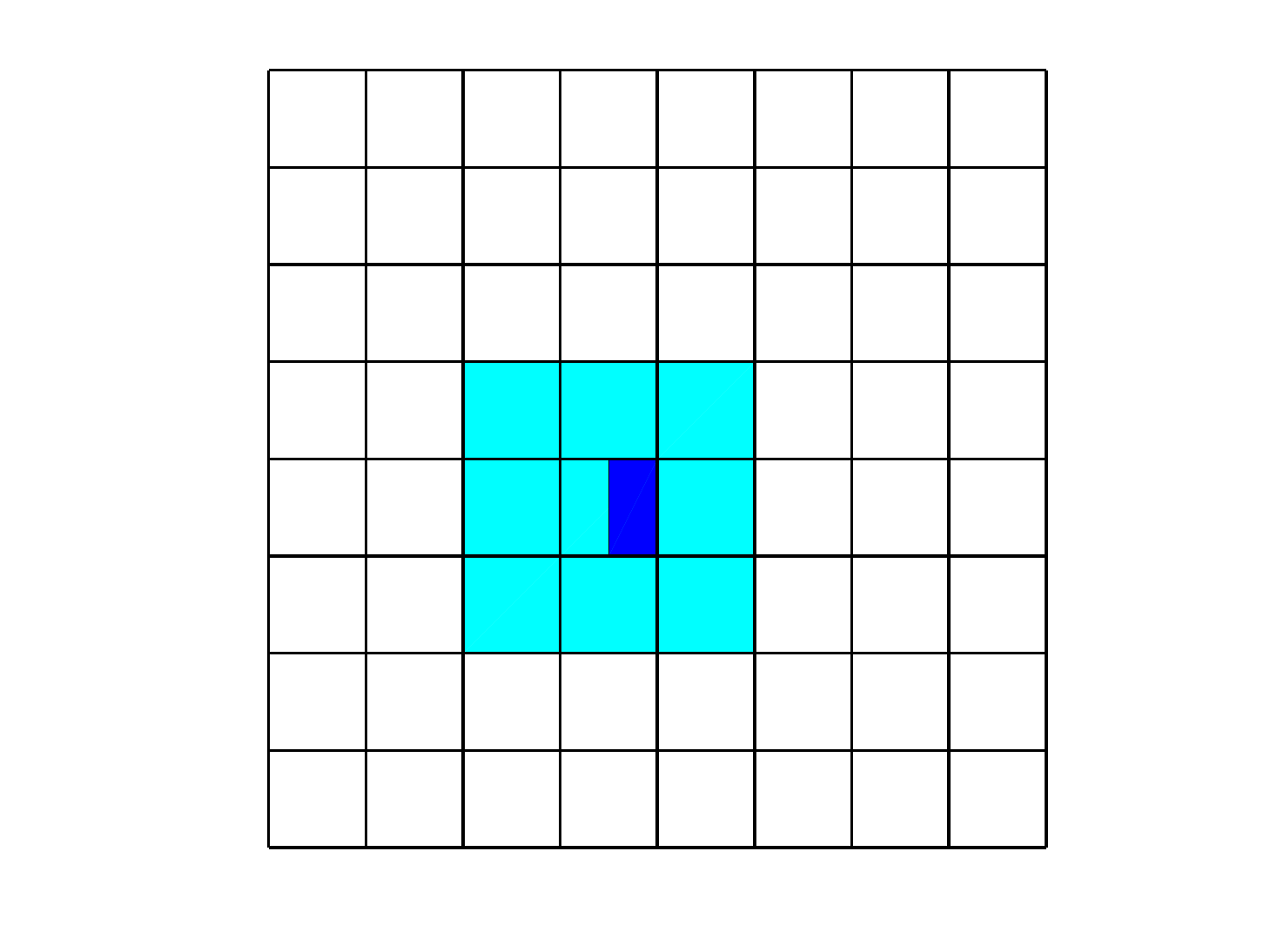}
\includegraphics[width=0.5\textwidth,trim=1cm 1cm 1cm 0cm, clip]{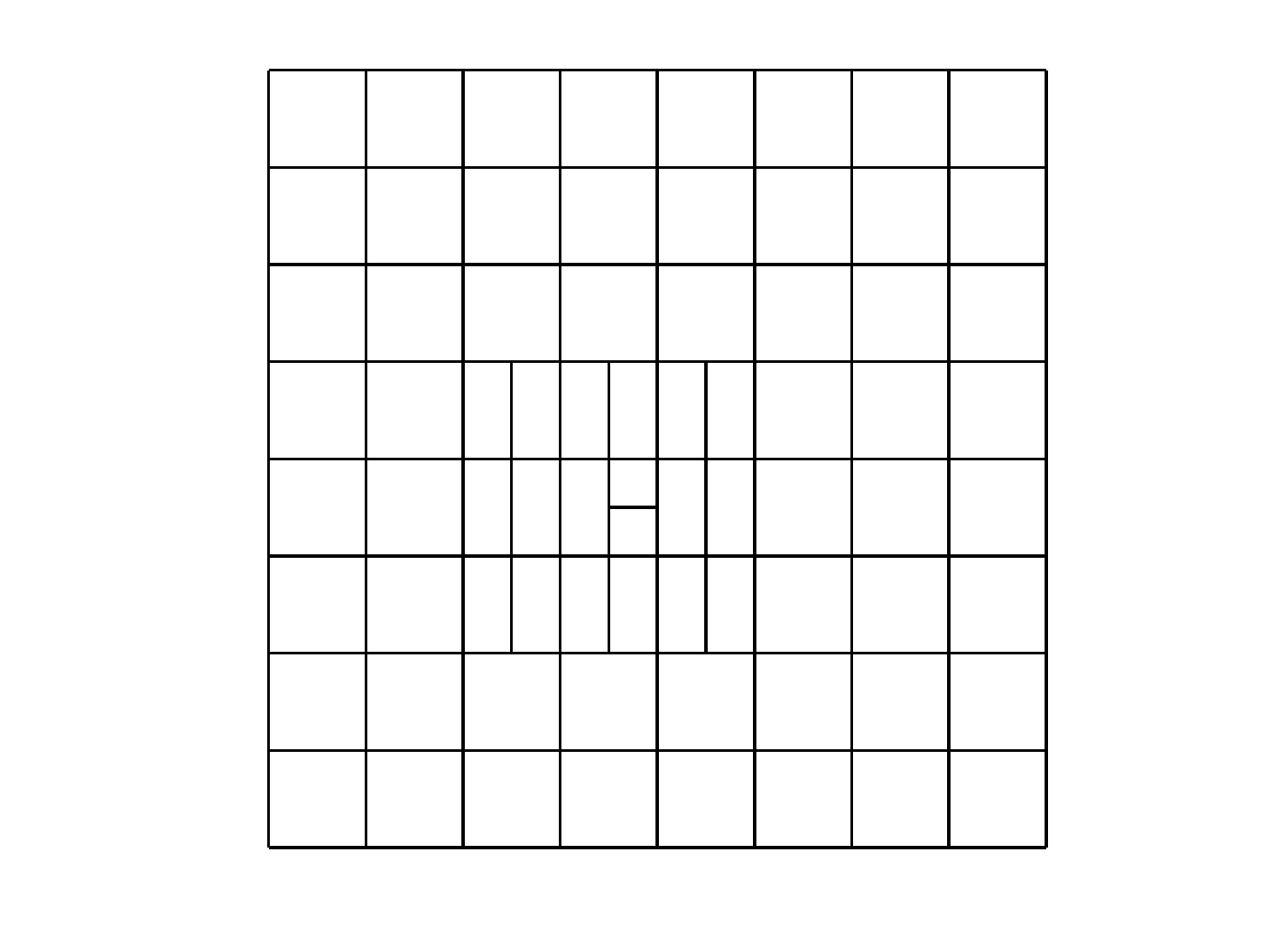}
\caption{$(3,3)$-neighborhood (in sky blue) of an element $\tau$ in blue (left). To obtain an admissible T-mesh, before bisecting $\tau$ it is necessary to bisect the elements in the neighborhood with lower generation (right)} \label{fig:neighborhood}
\end{figure}

Apart from the generations, that are necessary to provide the multilevel structure, we also need some definitions to adapt the refinement algorithm from \cite{Morgenstern_Peterseim} to the case of having an open knot vector. Given two points ${\bf x}, {\bf x'} \in \R^2$ we define the distance between them componentwise as the vector
\begin{equation*}
{\rm Dist} ({\bf x}, {\bf x'}) = ({\rm abs} \, (x_1 - x'_1), \, {\rm abs} \, (x_2 - x'_2)) \in \R^2.
\end{equation*}
Moreover, let us define, for a point in the index domain ${\bf x} = (x_1, x_2) \in [0, m_1] \times [0, m_2]$, its translated version ${\bf \tilde x} = (\tilde x_{1}, \tilde x_{2})$ given by
\begin{equation*}
\tilde x_{d} = \left \{
\begin{array}{ll}
p_d & \text{ if } x_{d} < p_d, \\
n_d & \text{ if } x_{d} > n_d, \\
x_{d} & \text{elsewhere}.
\end{array}
\right.
\end{equation*}
Given a bisection T-mesh ${\cal T}$ and an element $\tau \in {\cal T}$, we denote its middle point as ${\bf x}_\tau = (x_{1,\tau}, x_{2,\tau})$. We can define the distance of a point ${\bf x}$ to the element $\tau$, and the distance between two elements $\tau$ and $\tau'$ as
\begin{equation*}
{\rm Dist} ({\bf x}, \tau) := {\rm Dist} ({\bf \tilde x}, {\bf \tilde x}_\tau),
\qquad {\rm Dist} (\tau, \tau') := {\rm Dist} ({\bf \tilde x}_\tau, {\bf \tilde x}_{\tau'}),
\end{equation*}
respectively. Then, we define the $\bp$-neighborhood of $\tau$ as (see \cite[Definition~2.4]{Morgenstern_Peterseim} and Figure~\ref{fig:neighborhood})
\begin{equation*}
{\cal G}_\bp(\tau) : = \{ \tau' \in {\cal T} : {\rm Dist}(\tau, \tau') \le
D_\bp(g(\tau)) \},
\end{equation*}
where
\begin{equation*}
D_\bp(\ell) = \left \{
\begin{array}{ll}
2^{-\ell/2} \left( \lfloor p_1/2 \rfloor + 1/2, \; \lceil p_2/2 \rceil +
1/2
\right) & \text{ if $\ell$ is even} , \\
2^{-(\ell+1)/2} \left( \lceil p_1/2 \rceil + 1/2, \; 2\lfloor p_2/2 \rfloor +
1\right) & \text{ if $\ell$ is odd}.
\end{array}
\right.
\end{equation*}
We also need to define the set (see \cite[Corollary~2.15]{Morgenstern_Peterseim})
\begin{equation*}
U_\bp(\tau) = \{ {\bf x} \in [0,m_1] \times [0,m_2] : {\rm Dist}(\tau, {\bf x}) \le D_\bp(g(\tau)) \}.
\end{equation*}

Given a bisection T-mesh ${\cal T}$ and an element $\tau \in {\cal T}$, we say
that the bisection of $\tau$ is $\bp$-admissible or simply admissible, if all
$\tau' \in {\cal G}_\bp(\tau)$ satisfy $g(\tau') \ge g(\tau)$. Moreover,
we say that a bisection T-mesh ${\cal T}$ is $\bp$-admissible, if it can be
obtained as in \eqref{eq:T-by-bisection} with a sequence of admissible
bisections, see the example of Figure~\ref{fig:neighborhood}. It has been proved in \cite[Theorem~3.6]{Morgenstern_Peterseim} that admissible T-meshes are also analysis suitable, and therefore the dual basis and the projector of the previous section can be built.
\begin{remark}
The result in \cite{Morgenstern_Peterseim} does not take into account the repeated knots of the open knot vector. However, the same ideas apply using the definitions above, because the bisection of zero measure elements only adds new lines in the ``safe'' direction, without causing intersection of T-junction extensions. The use of the translated points in practice forces that a line arriving at a repeated knot, which is on the boundary of the parametric domain, should continue until the boundary of the index domain.
\end{remark}
Besides the analysis suitability condition, we need a local quasi-uniformity result, for which it is necessary to use the following auxiliary lemmas.
\begin{lemm} \label{lemma:bezier-tiled}
Let $\bA \in {\cal A}_\bp({\cal T})$ be an anchor of a $\bp$-admissible T-mesh ${\cal T}$. Then, any cell in the parametric tiled floor of $\bA$ contains at most two B\'ezier elements of ${\cal M^T}$.
\end{lemm}
{\it Proof.} Given a cell of the tiled floor, from \cite[Lemma~3.2]{BBSV2013} it does not contain any vertex of ${\cal T}$ in its interior, and any line of the extended mesh in its interior belongs to a T-junction extension. Since the mesh is analysis suitable, only vertical or horizontal lines can be found in its interior, but not both. Since a cell in the parametric tiled floor corresponds to a cell in the tiled floor, the result holds because we are refining by bisection alternating the refinement directions.
\hfill $\square$

\begin{lemm}\label{lemma:MP}
Given a $\bp$-admissible T-mesh ${\cal T}$ and an element $\tau \in {\cal T}$, for any $\tau' \in {\cal G}_\bp(\tau)$ it holds $g(\tau') \ge g(\tau) - 1$.
\end{lemm}
{\it Proof.} See \cite[Lemma~2.14]{Morgenstern_Peterseim}

\begin{lemm}\label{lemma:lqu}
Let $\bA \in {\cal A}_\bp({\cal T})$ be an anchor of a $\bp$-admissible T-mesh ${\cal T}$, associated to a function of generation $\ell$. Then the length of the cells of the parametric tiled floor of $\bA$ in each parametric direction is equal to
\begin{equation*}
L_1 =
\left \{
\begin{array}{cccl}
\frac{h_1}{2^{\ell/2}} &&& \text{ if } \ell \text{ is even}, \\
\frac{h_1}{2^{(\ell+1)/2}} &\text{ or }& \frac{h_1}{2^{(\ell-1)/2}}& \text{ if } \ell \text{ is odd},
\end{array}
\right.
\quad
L_2 =
\left \{
\begin{array}{cccl}
\frac{h_2}{2^{\ell/2}} &\text{ or }& \frac{h_2}{2^{\ell/2 - 1}} & \text{ if } \ell \text{ is even}, \\
\frac{h_2}{2^{(\ell-1)/2}} &&& \text{ if } \ell \text{ is odd}.
\end{array}
\right.
\end{equation*}

\end{lemm}
{\it Proof.} The lemma is proved by induction, similarly to \cite[Theorem~3.6]{Morgenstern_Peterseim}. The result clearly holds for the functions in the tensor-product space associated to the T-mesh ${\cal T}_0$. Assuming it is true for ${\cal T}_k$, we have to prove it for ${\cal T}_{k+1} = {\cal T}_k + b_{\tau_k}$. Since the bisections can be ordered by their generation, it holds that $\max_{\tau \in {\cal T}_k} g(\tau) = g(\tau_k) + 1$. We assume that $g(\tau_k)$ is even, the odd case is proved analogously exchanging the role of the vertical and horizontal directions.

Introducing ${\cal T}^u_\ell$ the Cartesian grid obtained after applying all possible bisections of generation $\ell$, using the same arguments as in \cite[Theorem~3.6]{Morgenstern_Peterseim}, we have the following result about the edges of the T-mesh in the $U_\bp(\tau_k)$ region
\begin{align}
\hsk({\cal T}_{k+1}) \cap U_\bp(\tau_k) = \hsk({\cal T}^u_{g(\tau_k)}) \cap U_\bp(\tau_k), \label{hsk-U}\\
\vsk({\cal T}^u_{g(\tau_k)}) \cap U_\bp(\tau_k) \subset \vsk({\cal T}_{k+1}) \cap U_\bp(\tau_k) \subset \vsk({\cal T}^u_{g(\tau_k)+1}) \cap U_\bp(\tau_k). \label{vsk-U}
\end{align}
We have to check that the result is true for the functions introduced or modified by the bisection $b_{\tau_k}$, that is, for functions in $\Phi_{k+1}$ of generation $\ell_{k+1} = g(\tau_k)+1$. Given one of these functions, from the length of the index vectors and \eqref{hsk-U}--\eqref{vsk-U} its associated anchor must intersect $U_\bp(\tau_k)$. Moreover, using \cite[Corollary~2.15]{Morgenstern_Peterseim} this anchor is associated to a geometrical entity (vertex, edge or element) that belongs to an element $\tau' \in {\cal G}_\bp(\tau_k)$. Using first Lemma~\ref{lemma:MP} and then the definition of $\bp$-admissible T-mesh, for any element $\tau'' \in {\cal G}_\bp(\tau')$ it holds $g(\tau'') \ge g(\tau')-1 \ge g(\tau_k) - 1$.

Since the horizontal index vector is collected from the vertical lines of elements in ${\cal G}_\bp(\tau')$, in the horizontal direction the intervals of the index vector have lengths at most $2^{-g(\tau_k)/2}$ and $2^{-(g(\tau_k)/2+1)}$, and recalling that $\ell_{k+1} = g(\tau_k) + 1$ is odd, we have the result for the horizontal direction.

For the vertical index vector, following an analogous reasoning we obtain that the intervals in the vertical direction also have lengths $2^{-g(\tau_k)/2}$ and $2^{-(g(\tau_k)/2+1)}$. Moreover, since the anchor of any function in ${\Phi}_{k+1}$ is horizontally aligned with the new edge, using~\eqref{hsk-U} the length of the vertical intervals is always $2^{-g(\tau_k)/2} = 2^{-(\ell_{k+1}-1)/2}$.

Finally, the result is proved passing from the index vectors to the local knot vectors, using the length of the intervals in the knot vectors $\Xi_d^\ell$ from Section~\ref{sec:T-bisection}.
\hfill $\square$

\begin{coroll}
In a $\bp$-admissible T-mesh, for any basis function ${B}_{\bA,\bp} \in \Phi_k$, and for any B\'ezier element $Q \subset {\rm supp}(B_{\bA,\bp})$, it holds that $h_Q \approx h_{\ell_k} \approx 2^{-\ell_k/2}$.
\end{coroll}

\begin{prop}\label{prop_hq_htildeq_hbarq}
For a $\bp$-admissible T-mesh ${\cal T}$, and for any B\'ezier element $Q \in  \cal M^T$ of the associated B\'ezier mesh, it holds that $h_Q \approx h_{\tilde Q} \approx h_{\bar Q}$.\\

\end{prop}
{\it Proof.} We only need to prove that $h_{\bar Q} \lec h_{\tilde Q} \lec h_Q$, because $Q \subset \tilde Q \subset \bar Q$. We start proving $h_{\tilde Q} \lec h_Q$. Since the T-mesh is analysis-suitable, in the element $Q$ there are at most $(p_1+1) (p_2+1)$ basis functions that do not vanish \cite[Proposition~7.6]{IGA-acta}, and for each of these functions the element $Q$ is contained in a cell of its parametric tiled floor. Moreover, from Lemma~\ref{lemma:bezier-tiled} this cell contains at most two elements of the B\'ezier mesh. Finally, each function contains at most $(p_1+1) (p_2+1)$ cells in its tiled floor, and the result of Lemma~\ref{lemma:lqu} states that the size of all these cells is comparable. Combining all these results we obtain $h_{\tilde Q} \lec h_Q$ with a constant that depends on the degree $\bp$.

Let us denote by $h_{1,Q}$ and $h_{2,Q}$ the length of $Q$ and in each parametric direction, and by $h_{1,\tilde Q}$ and $h_{2,\tilde Q}$ the length of $\tilde Q$, given by the difference in the first (second) parametric direction between the rightmost (uppermost) and leftmost (downmost) points in $\tilde Q$. From the previous results and the relative position of the B\'ezier element inside the tiled floor of each function (see the proof of \cite[Proposition~7.6]{IGA-acta}), we have that\[
(2p_d + 1) h_{d,Q} /2  \le h_{d,\tilde Q} \le (2 p_d + 1) 2 h_{d,Q},
\]
and since $\bar Q$ is the minimum rectangle that contains $\tilde Q$, it also  holds that $h_{\bar Q} \lec h_Q$.


\hfill $\square$

\begin{remark}
For simplicity we have assumed that the internal knots in the knot vectors $\Xi_d^\ell$ are equally spaced. The results of this section, and in particular Lemma~\ref{lemma:lqu}, can be extended to the case of local quasi-uniform knot vectors under Assumption~\ref{assumpt_quasiuniform}.
\end{remark}

\begin{figure}[ht]
\centering
\includegraphics[width=0.5\textwidth]{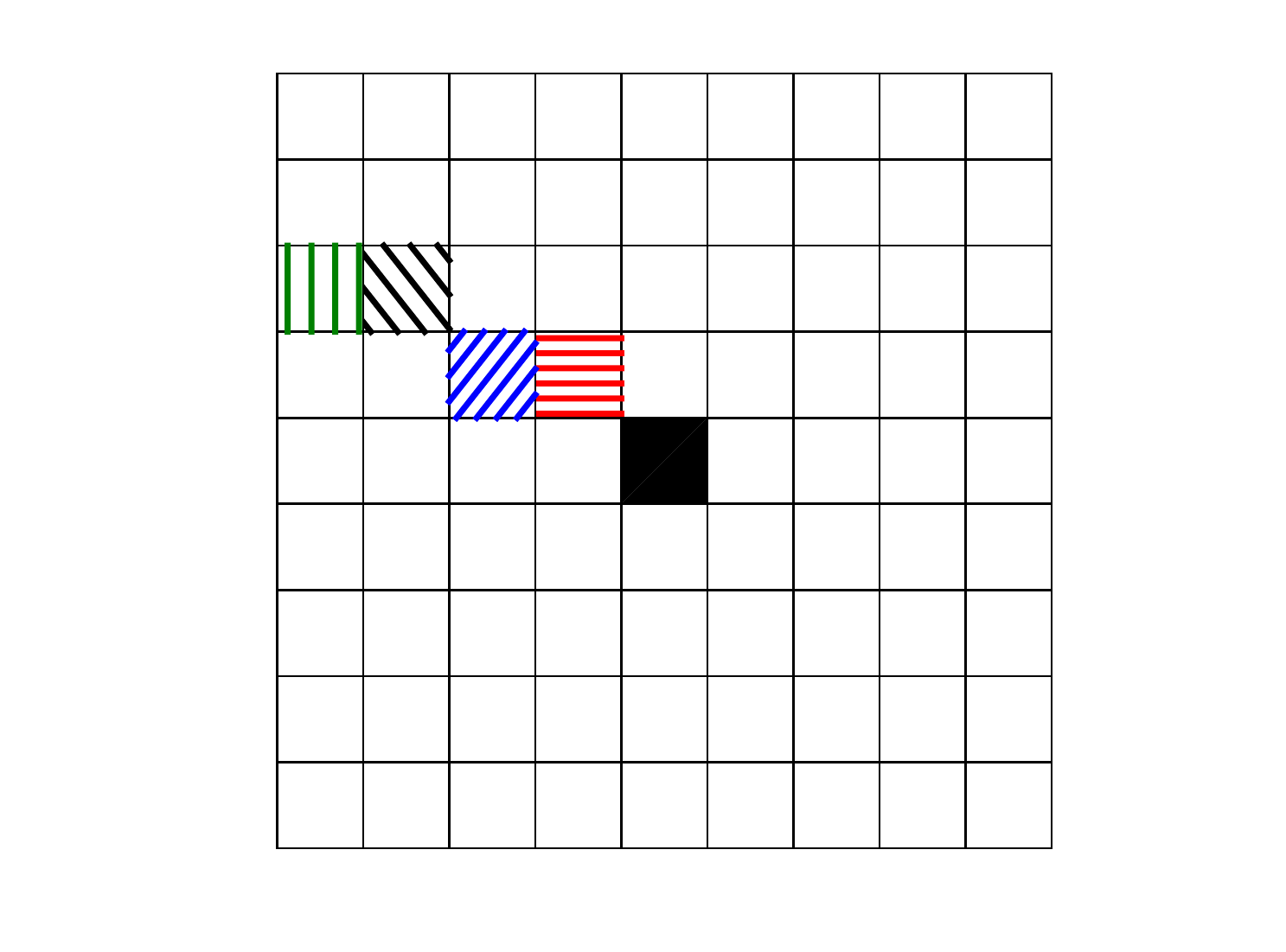}
\caption{For the chosen element $Q$ (in black), we represent the furthest element $\tau_k$ above and to the left such that, when $\tau_k$ is bisected, $Q$ is contained in $\omega_k$, for bilinear (red horizontally hatched), biquadratic (blue left diagonally hatched), bicubic (black right diagonally hatched) and biquartic (green vertically hatched).}
\label{fig:omega_k}
\end{figure}

\begin{prop} \label{prop:nonoverlapping}
For a ${\bf p}$-admissible T-mesh defined as in \eqref{eq:T-by-bisection}, and for all $u \in L^2(\Omega)$, and $\ell \ge 1$, we have that
\begin{equation}\label{nonoverlapping_consequence}
\sum_{k: \ell_k = \ell} \|u\|^2_{\omega_k} \lesssim \|u\|^2_{\Omega} \quad\mbox{and}\quad
\sum_{k: \ell_k = \ell} \|u\|^2_{\tilde{\omega}_k} \lesssim \|u\|^2_{\Omega},
\end{equation}
where $\omega_k$ and $\tilde \omega_k$ are as in \eqref{eq:omegak} and \eqref{eq:omegak-tilde}, respectively.
\end{prop}
{\it Proof.} The result follows from a careful counting of the number of times a B\'ezier element $Q$ (or its children, if it is bisected) can appear in the sets $\omega_k$ and $\tilde \omega_k$ of the same generation $\ell_k$. Studying the worst case scenario, when all the elements in the neighboring of $Q$ are bisected, we see from the graphical explanation in Figure~\ref{fig:omega_k} that the element is contained in at most $(2p_1+1)(2\lceil p_2/2 \rceil +1)$ (or $(2p_2+1)(2\lceil p_1/2 \rceil +1)$) sets $\omega_k$.

\begin{figure}
\centering
\includegraphics[width=0.32\textwidth]{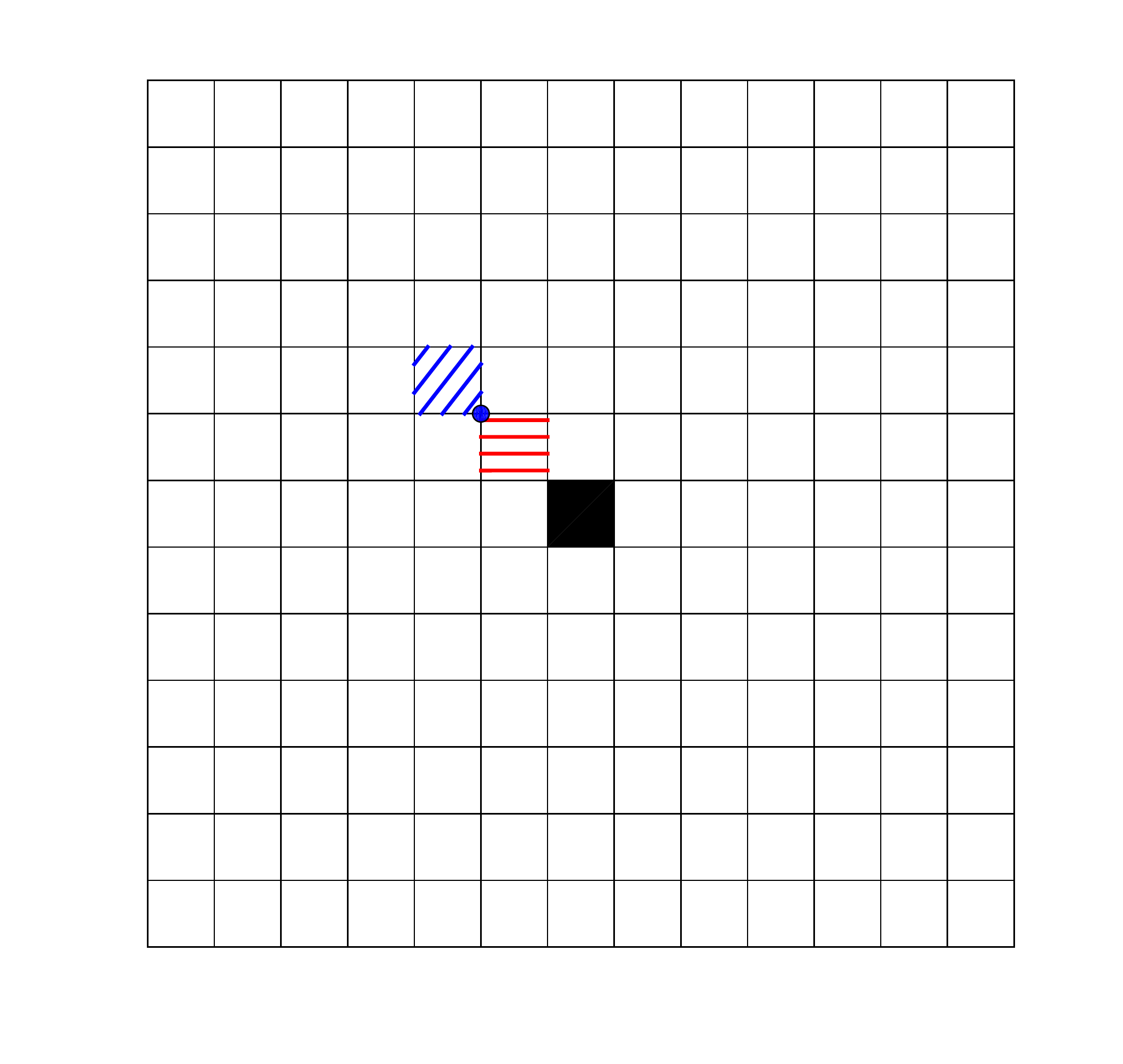}
\includegraphics[width=0.32\textwidth]{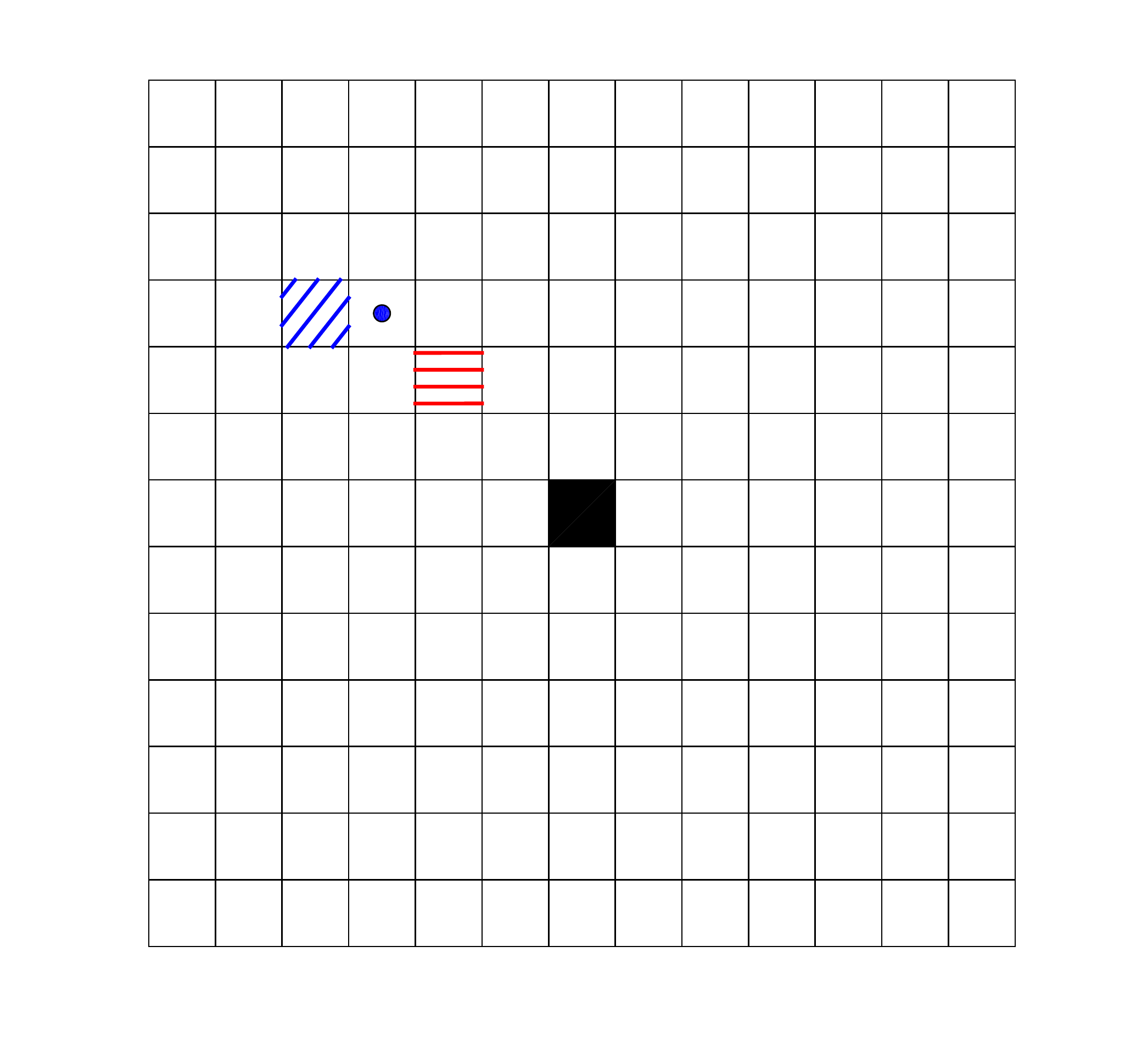}
\includegraphics[width=0.32\textwidth]{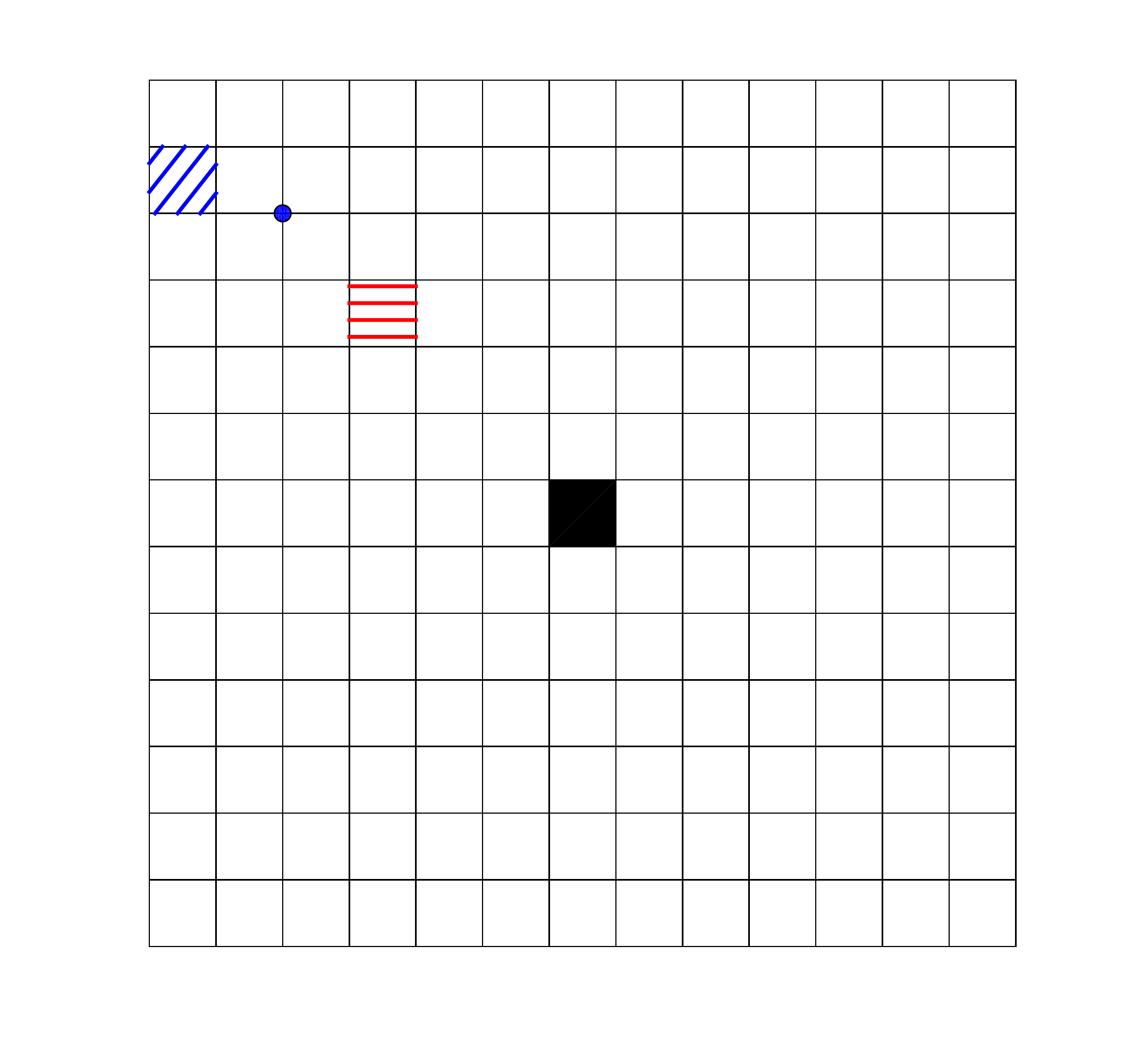}
\caption{For the chosen element $Q$ (in black), we represent the furthest element $Q'$ (horizontally hatched) such that there is a function that contains both elements in its support, and the furthest element (diagonally hatched) such that its bisection affects a function, anchored at the blue dot, that is in $\Phi_k$ and contains $Q'$ in its support. Thus, $Q' \subset \omega_k$ and $Q \subset \tilde \omega_k$. For the bilinear (left), biquadratic (middle) and bicubic (right) cases.}
\label{fig:omega_ktilde}
\end{figure}

A similar counting, that can be understood from Figure~\ref{fig:omega_ktilde}, gives that the element $Q$ is contained in at most $(4p_1+1)(4\lceil p_2/2 \rceil + 2 \lfloor p_2/2 \rfloor+1)$ (or $(4p_2+1)(4\lceil p_1/2 \rceil + 2 \lfloor p_1/2 \rfloor+1)$) sets $\tilde \omega_k$. The result holds with a constant that depends on the degree $\bp$.

We remark that in the figures we are mixing B\'ezier elements and elements in the T-mesh. In this worst case scenario, and staying away from the boundary, there is a one-to-one correspondence between both.\hfill $\square$

\begin{remark}
As already mentioned at the beginning of this section, the results can be extended to T-splines with lower continuity, that is, with repeated knots in the T-mesh. It would be necessary to adapt the definition of T-mesh by bisection in Section~\ref{sec:T-bisection} to take into account the multiplicity. The definition of analysis-suitable T-splines also applies to the case of reduced regularity, while the definition of {\bf p}-admissible T-meshes applies as well, taking into account that the neighborhood is defined in the T-mesh. We note that reducing the regularity reduces the support extension $\tilde Q$. Even if the proofs could be extended applying the same ideas, allowing for repeated knots would require a much more intricate notation, since one could not rely on the definition of the index sets ${\cal I}_d^\ell$ of each level.
\end{remark}

{\Bd
\begin{remark} \label{rem:constants}
The constants appearing in the proofs of the previous results, and in particular in Propositions~\ref{L2_stable_Pi_AST}, \ref{prop_hq_htildeq_hbarq} and~\ref{prop:nonoverlapping}, are independent of the mesh size but dependent on the degree. It is also important to note that the constants in \eqref{nonoverlapping_consequence} depend not only on the degree, but more precisely in the number of times that each single element appears in $\omega_k$ (and $\tilde \omega_k$), that is, in the overlaps between the differents $\omega_k$ of the same generation. This can have an important impact in the performance of the preconditioner, as we will see in Section~\ref{numerical_results}.
\end{remark}
}

\section{Space decomposition on AS T-meshes}\label{section5}
After having introduced all the preliminary results, we are now in a position to prove the main results of the paper, that is, a space decomposition that satisfies \Aone and \Atwo.

Let ${\cal T}$ be a $\bp$-admissible T-mesh obtained by successive bisections from the index T-mesh ${\cal T}_0$.
We give a decomposition of T-splines space $\mathcal{V}:=S(\mathcal{A}_{\bf p}(\mathcal{T}))$ using these successive bisections. By \eqref{subspace_Vk}, we have the following space decomposition:
\begin{equation}\label{space_decomposition}
\mathcal{V} =\sum_{k=0}^N \mathcal{V}_k.
\end{equation}
In each subspace $\mathcal{V}_k$, thanks to Lemma~\ref{lemma:lqu} and \cite{Bazilevs_Beirao_Cottrell_Hughes_Sangalli}, we have the following inverse inequality
\begin{equation}\label{inverse_inequality}
\|v_k\|_A^2 \lesssim h^{-2}_{\ell_k}  \|v_k\|_0^2,\quad \mbox{\ for\ all\ } v_k \in \mathcal{V}_k,
\end{equation}
where we recall that $h_{\ell_k}$ is the mesh size at level $\ell_k$, which is the level of functions in $\Phi_k$.


\subsection{Stable decomposition}

\begin{theo} [Space decomposition over AS T-meshes]\label{stabledecomp_ASTmeshes}
For any $v\in\mathcal{V}$, there exist $v_k\in {\cal V}_k,\ k=0,\ldots,N$ such that $v= \sum_{k=0}^N v_k$ and
\begin{equation}\label{ineq_spacedecomposition}
 \sum_{k=0}^N \|v_k\|^2_A \lesssim \|v\|_A^2.
\end{equation}

\end{theo}

{\it Proof.} First we consider an auxiliary decomposition over uniformly refined spaces of tensor-product B-splines. Let us recall that ${\cal T}_\ell^u$ indicates, for $0 \le \ell \le L$, the Cartesian mesh obtained after applying all possible bisections of generation $\ell$, $S_{\bf p}({\bf \Xi}^\ell)$ is the associated space of tensor-product B-splines of degree $\bp$, and ${\bf \Pi}_{\bf p, \Xi^\ell} : S_\bp({\bf \Xi}^L) \rightarrow S_\bp({\bf \Xi}^\ell)$ is the multivariate quasi-interpolant in \eqref{multi_quasiint}. We state (without proof) the following well-known stable decomposition for the space $S_\bp({\bf \Xi}^L) =\sum_{\ell=0}^L S_\bp({\bf \Xi}^\ell)$ (see \cite{BHKS13} for details)



\begin{lemm} [Stable decomposition for quasi-uniform meshes]\label{stabledecom_quasiuniform}
For any $\bar{v}\in S_\bp({\bf \Xi}^L)$, let $\bar{v}_\ell=({\bf \Pi}_{\bp, {\bf \Xi}^\ell}-{\bf \Pi}_{\bp, {\bf \Xi}^{\ell-1}})\bar{v}$
for $\ell=0,\ldots,L$, setting ${\bf \Pi}_{\bp, {\bf \Xi}^{-1}}:= 0$. Then $\bar{v}=\sum_{\ell=0}^L \bar{v}_\ell$ is a stable decomposition in the sense that
\begin{equation}
\sum_{\ell=0}^L h_\ell^{-2} \| \bar{v}_\ell \|^2 \lesssim |\bar{v}|_1^2.
\end{equation}
\end{lemm}

We now give a multilevel decomposition of $v$ using a sequence of projection operators
${\bf \Pi}_{\bf p}^{{\cal T}_k} : \mathcal{V} \rightarrow S(\mathcal{A}_{\bf p}(\mathcal{T}_k)),\ k=0,1,\ldots,N$, which are defined in \eqref{Projection_AST}.

First let $v = \sum_{\ell=0}^L \bar{v}_\ell$ with $\bar{v}_\ell :=({\bf \Pi}_{\bf p, \Xi^\ell}-{\bf \Pi}_{\bp, {\bf \Xi}^{\ell-1}})v \in S_\bp({\bf \Xi}^\ell)$. Then the slicing operator
${\bf \Pi}_{\bf p}^{{\cal T}_k}-{\bf \Pi}_{\bf p}^{{\cal T}_{k-1}}$ verifies an important property:
\[
v_k := ({\bf \Pi}_{\bf p}^{{\cal T}_k}-{\bf \Pi}_{\bf p}^{{\cal T}_{k-1}})v\ \in\ \mathcal{V}_k\qquad 0\le k\le N,
\]
with the convention that ${\bf \Pi}_{\bf p}^{{\cal T}_{-1}}:=0$ and thus $v=\sum_{k=0}^{N} v_k$.
If $\ell_k = g(\tau_{k-1}) + 1$ is the generation of the functions in $\Phi_k$, obtained after the bisection of the element $\tau_{k-1}$ in \eqref{eq:T-by-bisection}, then for $l \le \ell_k - 1 = g(\tau_{k-1})$ it holds that $({\bf \Pi}_{\bf p}^{{\cal T}_k}-{\bf \Pi}_{\bf p}^{{\cal T}_{k-1}})\bar{v}_l=0$, which is shown in details in Appendix \ref{appendixA1}.
Thus for $1\le k \le N$, we have
\[
v_k = ({\bf \Pi}_{\bf p}^{{\cal T}_k}-{\bf \Pi}_{\bf p}^{{\cal T}_{k-1}}) \sum_{l=\ell_k}^L \bar{v}_l.
\]
Also, Proposition~\ref{L2_stable_Pi_AST} shows that
\[
\|v_k\|^2_{\omega_k} \lesssim \big\|\sum_{l=\ell_k}^L \bar{v}_l \big\|_{\widetilde{\omega}_k}^2,\qquad 1\le k\le N.
\]
From the definition of the regions $\omega_k$, applying first the inequality above and then inequality~\eqref{nonoverlapping_consequence}, we infer that for $\ell\ge 1$
\[
\sum_{k: \ell_k = \ell} \| v_k \|^2 = \sum_{k: \ell_k = \ell} \|v_k\|_{\omega_k}^2
\lesssim
\sum_{k: \ell_k = \ell}  \big\| \sum_{l= \ell_k}^L \bar{v}_l \big\|_{\widetilde{\omega}_k}^2
\lesssim
\big\| \sum_{l=\ell}^L \bar{v}_l \big\|^2_{\Omega}
\lesssim
\sum_{l=\ell}^L \|\bar{v}_l\|^2_\Omega.
\]
Also we notice that when $\ell=0$ the space is tensor product, then $v_0 = \bar v_0$, and we trivially have
\[
\sum_{k: \ell_k = 0} \| v_k \|^2 = \|v_0\|^2_{\Omega}= \| \bar v_0 \|_\Omega^2 \le \sum_{l=0}^L \|\bar{v}_l\|_{\Omega}^2.
\]
Applying the discrete Hardy inequality \cite[Lemma~4.3]{CNX} with $s = \gamma^2 = 1/2$
to $a_\ell=\|\bar{v}_\ell\|^2$ and $b_\ell=\sum_{k: \ell_k=\ell}\|v_k\|^2$,
we obtain
\[
\sum_{\ell=0}^L h_\ell^{-2} \sum_{k: \ell_k=\ell} \|v_k\|^2 \lesssim \sum_{\ell=0}^L h_\ell^{-2} \|\bar{v}_\ell\|^2.
\]
Finally, Lemma \ref{stabledecom_quasiuniform} leads to
\[
\sum_{k=0}^N h_{\ell_k}^{-2} \|v_k\|^2 = \sum_{\ell=0}^L h_\ell^{-2} \sum_{k: \ell_k=\ell} \|v_k\|^2 \lesssim \sum_{\ell=0}^L h_\ell^{-2} \|\bar{v}_\ell\|^2
\lesssim |v|_1^2.
\]
Combining this inequality with \eqref{inverse_inequality} gives the desired estimate \eqref{ineq_spacedecomposition}. \hfill $\square$

\subsection{Strengthened Cauchy-Schwarz (SCS) inequality}
The proof of the SCS inequality for T-splines relies on the same result for tensor-product B-splines. Although multilevel methods for B-splines are now rather classical, we have not explicitly found the result for tensor-product splines in the literature. For this reason, we start this section proving the SCS inequality in the tensor-product case.

\begin{lemm} [SCS inequality for B-splines on globally quasi-uniform meshes]\label{SCS_quasiuniform}
For $u_i \in S_\bp({\bf \Xi}^i)$ and $u_j \in S_\bp({\bf \Xi}^j)$
with $j \ge i$, we have
\[
(u_i,u_j)_A \lesssim \gamma^{(j-i)/2} |u_i|_1 h_j^{-1}\|u_j\|_0,
\]
where $\gamma < 1$ is a constant such that $h_i\approx \gamma^{i}$.
\end{lemm}

{\it Proof.}
Recalling that $S_\bp({\bf \Xi}^i)=\mbox{span}\{{B}_{\bf i,p}({\bzeta}),~ {\bf i} \in {\bf I}^i\}$ and $S_\bp({\bf \Xi}^j)=\mbox{span}\{{B}_{\bf j,p}({\bzeta}),~ {\bf j} \in {\bf I}^j\}$, any $u_j\in S_\bp({\bf \Xi}^j)$ can be written as
\[
u_j = \sum_{{\bf j} \in {\bf I}^j} c_{\bf j,p}{B}_{\bf j,p},
\]
which can also be classified as follows: for each B\'ezier element $Q^i \in {\cal M}_i^{\cal B}$ we define the sets of indices
\begin{align*}
&{\bf I}^j_{\rm out} := \{ {\bf j} \in {\bf I}^j : {\rm supp}({B}_{\bf j,p})\cap Q^i=\emptyset \}, \\
&{\bf I}^j_{\rm in} := \{ {\bf j} \in {\bf I}^j : {\rm int (supp}({B}_{\bf j,p})) \subset Q^i \}, \\
&{\bf I}^j_{\rm B} := \{ {\bf j} \in {\bf I}^j : {\rm int (supp}({B}_{\bf j,p}))\cap \partial Q^i \ne \emptyset \},
\end{align*}
of basis functions with support completely outside of $Q^i$, completely contained in $Q^i$, or with just part of it in $Q^i$, where we recall that the element $Q^i$ is open. Then, we have
\[
u_j = u_{\rm out}+u_{\rm in}+u_{\rm B} = \sum_{{\bf j} \in {\bf I}^j_{\rm out}} c_{\bf j,p}{B}_{\bf j,p} + \sum_{{\bf j} \in {\bf I}^j_{\rm in}} c_{\bf j,p}{B}_{\bf j,p} + \sum_{{\bf j} \in {\bf I}^j_{\rm B}} c_{\bf j,p}{B}_{\bf j,p}.
\]

Denoting $a_{Q^i}(u,v):=\int_{Q^i} \nabla u\cdot \nabla v~dx$, we clearly have
$a_{Q^i}(u_i,u_{\rm out})=0$, and thus
\[
a_{Q^i}(u_i,u_j)=a_{Q^i}(u_i,u_{\rm in})+a_{Q^i} (u_i,u_{\rm B}).
\]
Define $\Gamma := \bigcup_{{\bf j} \in {\bf I}^j_{\rm B}} {\rm supp} (B_{\bf j, p}) \cap Q^i$, which is the union of smaller elements $Q^j \in {\cal M}_j^{\cal B}$ contained in $Q^i$ and in the support of a function that is not completely contained in $Q^i$. It is easy to see that $| \Gamma| \approx h_i h_j$, with the implicit constant depending on the degree.

First we estimate the term corresponding to $u_{\rm B}$. Since $\mbox{supp}(u_{\rm B})\cap Q^i=\Gamma$, we have
\[
|a_{Q^i}(u_i,u_{\rm B})|=|a_{\Gamma}(u_i,u_{\rm B})| \le |u_i|_{H^1(\Gamma)} |u_{\rm B}|_{H^1(\Gamma)}.
\]
Using that $\Gamma \subset Q^i$ with $|\Gamma| \approx h_i h_j$, the inverse inequality \cite[Lemma~4.5.3]{Brenner_Scott_1994}, and finally that $h_j \approx \gamma^j$, we obtain
\begin{eqnarray*}
|u_i|_{H^1(\Gamma)} & \le & |\Gamma|^{1/2} |u_i|_{W^{1,\infty}(\Gamma)} \le |\Gamma|^{1/2} |u_i|_{W^{1,\infty}(Q^i)}\\
& \lesssim &  |\Gamma|^{1/2} h_i^{-1} |u_i|_{H^1(Q^i)} \lesssim h_j^{1/2} h_i^{1/2} h_i^{-1}  |u_i|_{H^1(Q^i)}\\
& \lesssim & \gamma^{(j-i)/2} |u_i|_{H^1(Q^i)}.
\end{eqnarray*}
The inverse inequality on globally quasi-uniform meshes, followed by \cite[Proposition~5.1 and Corollary~5.1]{BePa_BDDC} shows that
\[
|u_{\rm B}|_{H^1(\Gamma)}= \sum_{\substack{Q^j \subset \Gamma \\ Q^j \in {\cal M}^{\cal B}_j}} |u_{\rm B}|_{H^1(Q^j)}
\lesssim h_j^{-1} \sum_{\substack{Q^j \subset \Gamma \\ Q^j \in {\cal M}^{\cal B}_j}} \|u_{\rm B}\|_{L^2(Q^j)}
\lesssim h_j^{-1} \|u_{\rm B}\|_{L^2(\Gamma)} \lesssim h_j^{-1} \|u_j\|_{L^2(Q^i)},
\]
where a detailed proof of the last inequality can be found in Appendix \ref{appendixA2}.
Combining the above two inequalities leads to
\begin{equation}\label{ineq_uB}
|a_{Q^i}(u_i, u_{\rm B})| \lesssim \gamma^{(j-i)/2} |u_i|_{H^1(Q^i)} h^{-1}_j \|u_j\|_{L^2(Q^i)}.
\end{equation}

Now, for the estimate of the term corresponding to $u_{\rm in}$ we have
\begin{equation*}
|a_{Q^i}(u_i, u_{\rm in})| \le  |(-\Delta u_i, u_{\rm in})_{L^2(Q^i)}|
 \le  \|\Delta u_i\|_{L^2(Q^i)} \|u_{\rm in}\|_{L^2(Q^i)}.
\end{equation*}
It follows from the inverse inequality that
\[
\|\Delta u_i\|_{L^2(Q^i)} \lesssim |u_i|_{H^2(Q^i)} \lesssim h^{-1}_i |u_i|_{H^1(Q^i)}. 
\]
Noting that $\|u_{\rm in}\|_{L^2(Q^i)} \lesssim \|u_j\|_{L^2(Q^i)}$ that can be shown in a similar way to Appendix \ref{appendixA2}, we get
\begin{equation}\label{ineq_uIN}
\begin{array}{lll}
|a_{Q^i}(u_i,u_{\rm in})| & \lesssim & h^{-1}_i |u_i|_{H^1(Q^i)} \|u_j\|_{L^2(Q^i)} \\
& \lesssim & h_j h_i^{-1} h_j^{-1} |u_i|_{H^1(Q^i)} \|u_j\|_{L^2(Q^i)}\\
& \lesssim & \gamma^{j-i} |u_i|_{H^1(Q^i)} h_j^{-1}\|u_j\|_{L^2(Q^i)}.
\end{array}
\end{equation}
Combining the two estimates \eqref{ineq_uB} and \eqref{ineq_uIN}, together with $0 < \gamma < 1$, yields
\[
|a_{Q^i}(u_i,u_j)| \lesssim \gamma^{(j-i)/2} |u_i|_{H^1(Q^i)} h^{-1}_j \|u_j\|_{L^2(Q^i)}.
\]
Summing over $Q^i$ and using the \Bd Cauchy-Schwarz \Bk inequality, we obtain
\[
a(u_i,u_j) \lesssim \gamma^{(j-i)/2} |u_i|_1 h^{-1}_j \|u_j\|_0.
\]
which ends the proof.
\hfill $\square$

\begin{theo} [SCS inequality for AS T-meshes]\label{SCS_AST}
For any $u_i,v_i\in \mathcal{V}_i,\ 0\le i \le N$, we have
\[
\left| \sum_{i=0}^N \sum_{j=i+1}^N (u_i,v_j)_A \right|
\lesssim \left( \sum_{i=0}^N \|u_i\|_A^2 \right)^{1/2} \left( \sum_{i=0}^N \|v_i\|_A^2 \right)^{1/2}.
\]
\end{theo}
{\it Proof.} The proof is similar to the proof of Theorem 4.6 in \cite{CNX}, and is divided into four steps.

1. For a fixed $i\in\{1,2,\ldots,N\}$, we denote
\[
n(i):=\{j>i:~ \omega_j \cap \omega_i \neq \emptyset \} \quad
\mbox{\ and\ } w_{\ell}^i : = \sum_{\substack{j\in n(i) \\ j:\ell_j=\ell}} v_j.
\]
Observe that $w_{\ell}^i \in S_\bp({\bf \Xi}^{\ell})$ and $\ell=\ell_j \ge \ell_i$ by the assumption that if $k > k'$, then $\ell_k \ge \ell_{k'}$, while  $u_i\in S_\bp({\bf \Xi}^{\ell_i})$.

For any B\'ezier element $Q \subset {\omega}_i$, we apply Lemma \ref{SCS_quasiuniform} over $Q$ to $u_i$ and $w_{\ell}^i$ to obtain
\[
(u_i,w_{\ell}^i)_{A,Q} \lesssim \gamma^{(\ell-\ell_i)/2}\|u_i\|_{A,Q} \, h^{-1}_{\ell} \|w_{\ell}^i\|_{0,Q}.
\]
Then
\begin{eqnarray*}
(u_i,w_{\ell}^i)_{A,\omega_i} = \sum_{Q\subset\omega_i} (u_i,w_{\ell}^i)_{A,Q}
&\lesssim& \gamma^{(\ell-\ell_i)/2} \sum_{Q\subset \omega_i} \|u_i\|_{A,Q} \, h_{\ell}^{-1} \|w_{\ell}^i\|_{0,Q} \\
&\lesssim& \gamma^{(\ell-\ell_i)/2} \|u_i\|_{A,\omega_i} \, h_{\ell}^{-1} \big( \sum_{Q \subset \omega_i} \|w_{\ell}^i\|_{0,Q}^2 \big)^{1/2}.
\end{eqnarray*}

Moreover, from the definition of $w_\ell^i$
\[
\|w_{\ell}^i\|_{0,Q}^2 \le \sum_{\substack{j\in n(i) \\ j:\ell_j=\ell}} \|v_j\|_{0,Q}^2 \le \sum_{j:\ell_j=\ell} \|v_j\|_{0,Q}^2.
\]
Thus we obtain
\[
(u_i,w_{\ell}^i)_{A,\omega_i} \lesssim \gamma^{(\ell-\ell_i)/2} \|u_i\|_{A,\omega_i} h_{\ell}^{-1} \big( \sum_{j:\ell_j=\ell} \|v_j\|^2_{0,\omega_i} \big)^{1/2}.
\]

2. Fix $u_i$ and consider
\begin{eqnarray*}
\big| (u_i,\sum^N_{j=i+1} v_j)_A \big| &=& \big| (u_i, \sum_{j\in n(i)} v_j)_{A,\omega_i} \big|\\
&=& \big|  (u_i, \sum^L_{\ell=\ell_i}\sum_{\substack{j\in n(i)\\ j:\ell_j=\ell}} v_j)_{A,\omega_i} \big| \le
\sum^L_{\ell=\ell_i} \big| (u_i,w_{\ell}^i)_{A,\omega_i} \big|.
\end{eqnarray*}
Thus we get
\[
\big| (u_i,\sum^N_{j=i+1} v_j)_A \big| \lesssim
\sum^L_{\ell=\ell_i} \gamma^{(\ell-\ell_i)/2} \|u_i\|_{A,\omega_i} h_{\ell}^{-1} \big( \sum_{j:\ell_j=\ell} \|v_j\|^2_{0,\omega_i} \big)^{1/2}.
\]

3. We sum over $i$ with fixed generation $\ell_i=k$:
\begin{eqnarray*}
\sum_{i:\ell_i=k} \big| (u_i,\sum^N_{j=i+1} v_j)_A \big|
&\lesssim&
\sum_{i:\ell_i=k} \big( \sum^L_{\ell=\ell_i} \gamma^{(\ell-\ell_i)/2} \|u_i\|_{A,\omega_i} h_{\ell}^{-1} \big( \sum_{j:\ell_j=\ell} \|v_j\|^2_{0,\omega_i} \big)^{1/2} \big)\\
&\lesssim& \sum^L_{\ell=k} \gamma^{(\ell-k)/2}\Big( \sum_{i:\ell_i=k}\big[ \|u_i\|_{A,\omega_i} (h_{\ell}^{-2} \sum_{j:\ell_j=\ell} \|v_j\|^2_{0,\omega_i})^{1/2}
\big]\Big)\\
&\lesssim&
\sum^L_{\ell=k} \gamma^{(\ell-k)/2}  \Big( \sum_{i:\ell_i=k} \|u_i\|_{A,\omega_i}^2 \Big)^{1/2}
\Big( h_\ell^{-2} \sum_{i:\ell_i=k} \sum_{j:\ell_j=\ell} \|v_j\|^2_{0,\omega_i} \Big)^{1/2} \\
&\lesssim&
\sum^L_{\ell=k} \gamma^{(\ell-k)/2}  \Big( \sum_{i:\ell_i=k} \|u_i\|_{A,\omega_i}^2 \Big)^{1/2}
\Big( h_\ell^{-2} \sum_{j:\ell_j=\ell} \|v_j\|^2_{0} \Big)^{1/2}
\end{eqnarray*}
where in the third inequality we used the usual Cauchy-Schwarz inequality and in the last inequality we used Proposition~\ref{prop:nonoverlapping}.

4. Finally we sum over all the generations $0\le k \le L$ to get
\begin{eqnarray*}
\sum_{k=0}^L \sum_{i:\ell_i=k} \big| (u_i,\sum^N_{j=i+1} v_j)_A \big|
&\lesssim&
\sum_{k=0}^L \left( \sum^L_{\ell=k} \gamma^{(\ell-k)/2}  \Big( \sum_{i:\ell_i=k} \|u_i\|_{A,\omega_i}^2 \Big)^{1/2}
\Big( h_\ell^{-2} \sum_{j:\ell_j=\ell} \|v_j\|^2_{0} \Big)^{1/2} \right)\\
&\lesssim&
\Big( \sum_{k=0}^L \sum_{i:\ell_i=k} \|u_i\|_{A,\omega_i}^2 \Big)^{1/2}
\Big( \sum_{\ell=0}^L h_\ell^{-2} \sum_{j:\ell_j=\ell} \|v_j\|^2_0 \Big)^{1/2},
\end{eqnarray*}
where in the last inequality we used the inequality (see \cite[Lemma~5]{XCN})
\[
\sum_{i,j=1}^n \gamma^{|i-j|} x_i y_j \lesssim \frac{2}{1-\gamma} \big(\sum_{i=1}^n x_i^2 \big)^{1/2} \big(\sum_{j=1}^n y_j^2 \big)^{1/2}, \quad \forall (x_i)_{i=1}^n, (y_j)_{j=1}^n \in {\mathbb R}^n.
\]
Since $\sum^L_{k=0}\sum_{i:\ell_i=k}=\sum^N_{i=0}$,
using a scaled Poincar\'e inequality, namely,
\[
h^{-2}_{\ell_k} \|v_k\|^2_0 \lesssim  \|v_k\|_A^2, \quad \mbox{\ for\ all\ } v_k \in \mathcal{V}_k,
\]
we have
\[
\big| \sum_{i=0}^N \sum_{j=i+1}^N (u_i,v_j)_A \big| \le
\sum_{i=0}^N \big| (u_i,\sum_{j=i+1}^N v_j)_A  \big| \lesssim
\big( \sum_{i=0}^N \|u_i\|^2_A \big)^{1/2} \big( \sum_{j=0}^N \|v_j\|^2_A  \big)^{1/2}.
\]
\hfill $\square$

\section{BPX preconditioners on locally quasi-uniform AS T-meshes}
After the proof of the stable decomposition \Aone and the SCS inequality \Atwo, we can apply the parallel subspace correction method. We introduce two different space decompositions: the first one, that we called \emph{micro decomposition}, is based on the same subspaces in the previous two sections; the second one, that we denote \emph{macro decomposition}, is based on a smaller number of subspaces, that collect all the bisections of the same generation at once.

\subsection{Micro decomposition}\label{sec:micro}
We apply parallel subspace correction methods to the space decomposition
\begin{equation}\label{space_decomposition2}
\mathcal{V} =\sum_{k=0}^N \mathcal{V}_k,
\end{equation}
where the subspaces ${\cal V}_k$ are defined in \eqref{subspace_Vk}, and thus obtain BPX preconditioners on locally quasi-uniform AS T-meshes.

\begin{theo} If the preconditioner $B$ of \eqref{def_B} is based on the space
decomposition \eqref{space_decomposition2} and the SPD smoothers satisfying \eqref{assumption_smoother}, then we have
\[
\kappa(BA)\lesssim 1.
\]
\end{theo}

{\it Proof.}  It follows immediately from Theorems \ref{stabledecomp_ASTmeshes}, \ref{SCS_AST} and \ref{BPX_preconditioning}. \hfill $\square$

\begin{remark} We observe that we can use standard smoothers such as Jacobi and symmetric Gauss-Seidel iterations satisfying \eqref{assumption_smoother}. Recently, Hofreither et al. in \cite{Hofreither2016} proposed a new smoother based on the mass matrix and a boundary correction, {\Bd that was later improved in \cite{Hofreither-Takacs}. The results in those papers show that on a uniform mesh, the multigrid with the new smoother is more robust, in the sense that convergence is independent of both the mesh size and the spline degree.} The extension of this smoother to T-splines is one of our future research topics. We notice however that this would eliminate the dependency on the degree for $K_3$ and $K_4$ in \eqref{eq:cond}, but not for $K_1$ and $K_2$.
\end{remark}

\subsection{Macro decomposition}\label{sectionMacro}
As an alternative to the previous decomposition, we introduce a \emph{macro} space decomposition where each level contains all the elements of the same corresponding generation. In order to define the macro space decomposition, to each generation $\ell$ we associate a subspace ${\mathcal W}_\ell$ that contains the functions added or modified after inserting all the lines of generation $\ell$. For instance, in the T-meshes of Figure~\ref{fig3:meshes}, ${\mathcal W}_0$ would be the tensor-product space in the \emph{black} mesh, and the subspaces ${\mathcal W}_1$, ${\mathcal W}_2$ and ${\mathcal W}_3$ would consist of the T-splines added or modified after the bisection sequence for adding all the \emph{red} vertical edges, all the \emph{blue} horizontal edges, and all the \emph{pink} vertical edges, respectively.

To be more precise, starting from the Cartesian mesh ${\cal T}_0$, we define for $1 \le \ell \le L$ the T-mesh
\begin{equation*}
{\cal \overline T}_\ell := {\cal T}_{k(\ell)}, \qquad \text{ with } \quad k(\ell) := \max \{ k : \ell_k=\ell \mbox{\ with\ } \ell_k := g(b_{\tau_{k-1}})\},
\end{equation*}
which is the finest T-mesh of generation $\ell$. Similarly to the sets $\Phi_k$ in Section~\ref{sec:ASTS-bisection}, we can now define the sets
\[
{\overline \Phi}_\ell := \{{B}_{\bf A,p}: {\bf A}\in \mathcal{A}_{\bf p}(\mathcal{\overline T}_{\ell})\} \backslash \{{B}_{\bf A,p}: {\bf A}\in \mathcal{A}_{\bf p}(\mathcal{\overline T}_{\ell-1})\},
\]
and the subspace
\[
{\cal W}_\ell := {\rm span} \; {\overline \Phi}_\ell,
\]
noting that this subspace contains all the subspaces of the micro decomposition of the corresponding generation, in fact, ${\cal W}_\ell = \sum_{k:\ell_k = \ell} {\cal V}_k$.

Defining those ${\mathcal W}_\ell$ gives the decomposition
\begin{equation}
\mathcal{V} =\sum_{\ell=0}^{L} \mathcal{W}_\ell.\label{macrodecomp}
\end{equation}
We can use the general framework presented in Sections~\ref{section4}-\ref{section5} to obtain the uniformly bounded condition number of the BPX preconditioner  on the macro space decomposition \eqref{macrodecomp}.

\begin{remark} It is worth to mention that since our macro structure is similar to the locally refined grids proposed in \cite{BPWX,BPX1,Xu_MGintro}, BPX optimality on the space decomposition \eqref{macrodecomp} can also be obtained in an analogous way to the proofs in those references.
\end{remark}

\subsection{Application to IGA} \label{sec:IGA}
All the theoretical results presented in this paper are proved in the unit domain ${\Omega} = (0,1)^2$. The results can be easily extended to the IGA setting, where the physical domain $\Omega$ is defined as the image of the unit parametric domain through a parametrization ${\bf F}$, that is, $\Omega = {\bf F}((0,1)^2)$. The parametrization ${\bf F}$ can be defined as a T-spline on a coarse mesh, and both ${\bf F}$ and its inverse should be regular, see \cite{BHKS13} for details.

\section{Numerical results}
\label{numerical_results}
We have performed some numerical tests to support our theoretical results with numerical evidence. We solve the model problem \eqref{model_prob} with a suitable $f$ that makes it necessary to refine towards a corner, noticing however that the condition number does not depend on $f$. Since we lack a true implementation of T-splines, we have run the tests using the Octave/Matlab software GeoPDEs \cite{dFRV11,GEOPDES-NEW} for some particular T-meshes, defining a tensor-product space for each level, and then collecting the active functions of different levels. Our implementation is clearly inefficient for T-splines, hence we do not present computational times.

\subsection{Macro tests} \label{sec:macro-tests}
In the first two tests we follow the \emph{macro} space decomposition in Section~\ref{sectionMacro}, that is, the subspace ${\cal W}_k$ contains the functions that have been added or modified after inserting all the lines of generation $k$.

\subsubsection{Square domain}
\begin{figure}[htp!]
\begin{subfigure}[T-mesh, $p=2$]{
\includegraphics[width=0.48\textwidth]{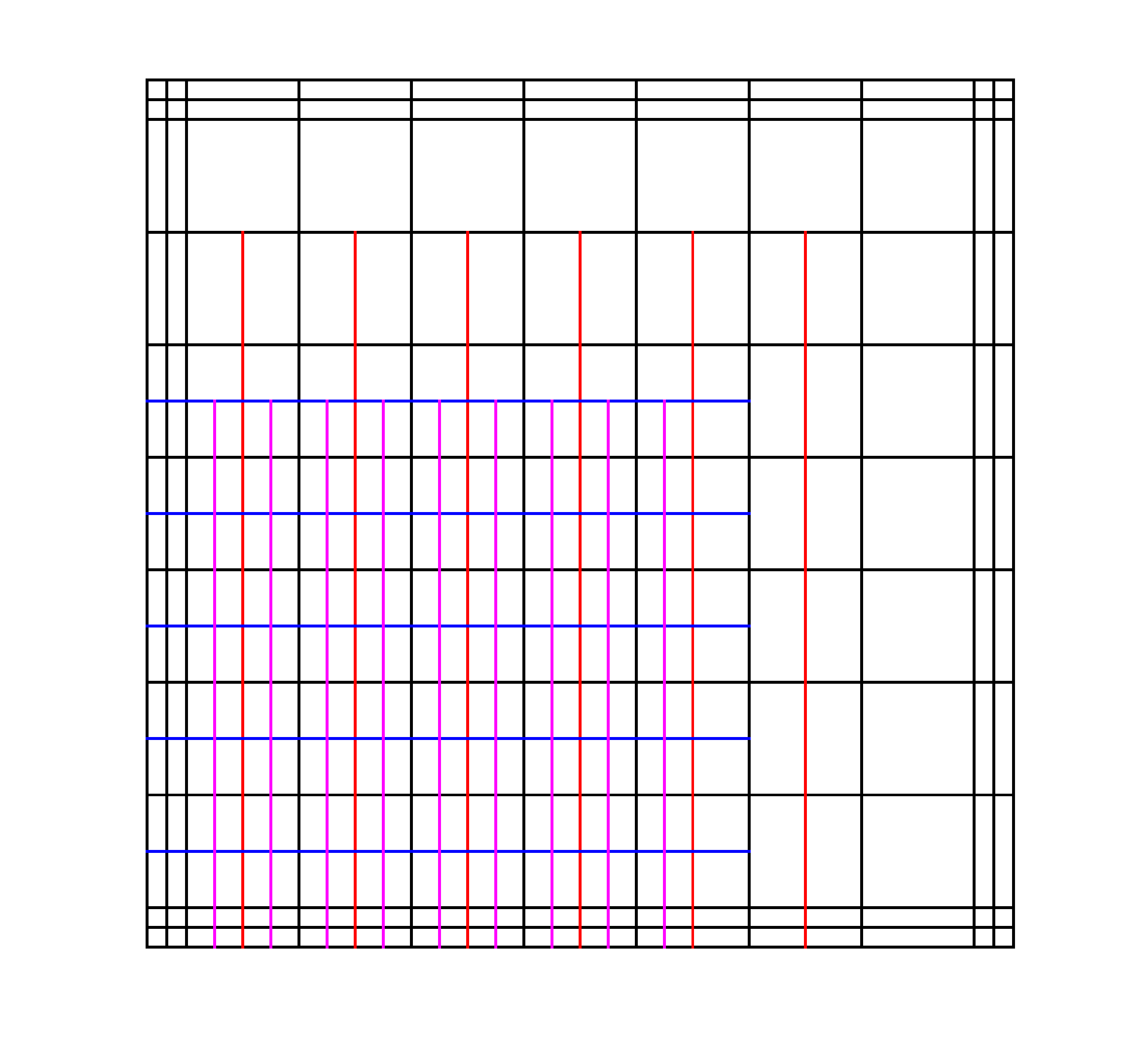}
}
\end{subfigure}
\begin{subfigure}[B\'ezier mesh, $p=2$]{
\includegraphics[width=0.48\textwidth]{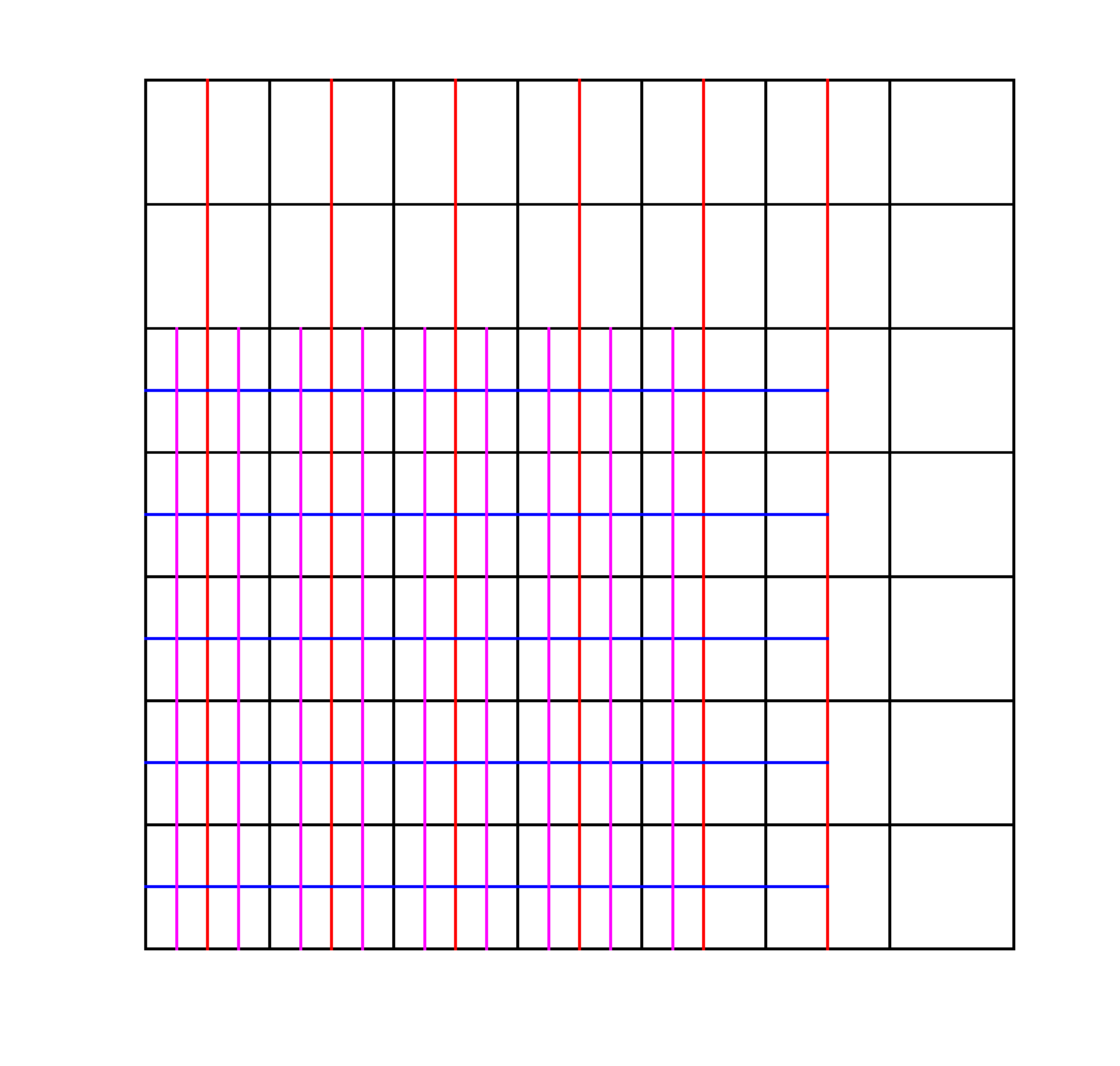}
}
\end{subfigure}
\begin{subfigure}[T-mesh, $p=3$]{\label{fig:Tmesh-square-p3}
\includegraphics[width=0.48\textwidth]{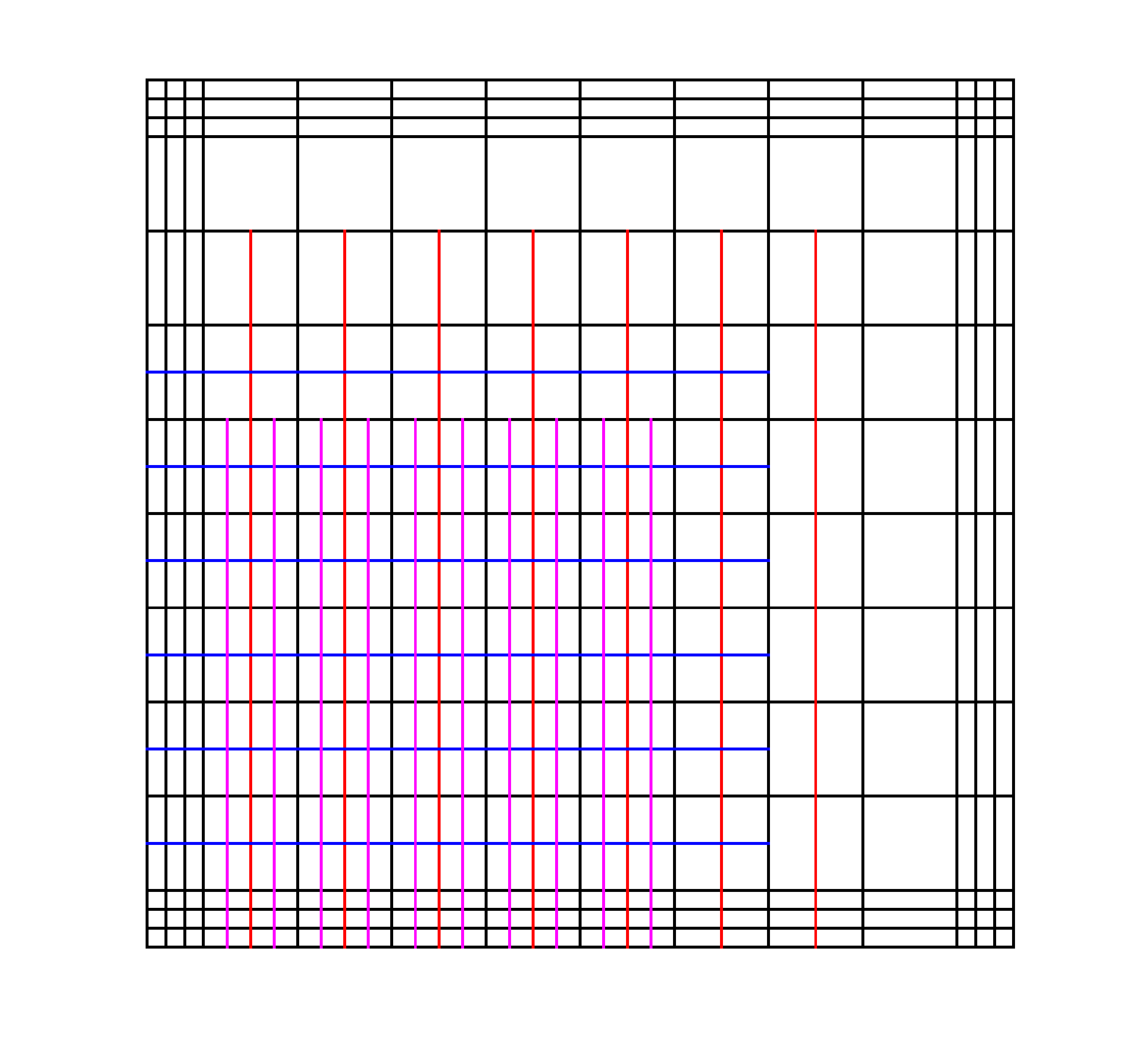}
}
\end{subfigure}
\begin{subfigure}[B\'ezier mesh, $p=3$]{
\includegraphics[width=0.48\textwidth]{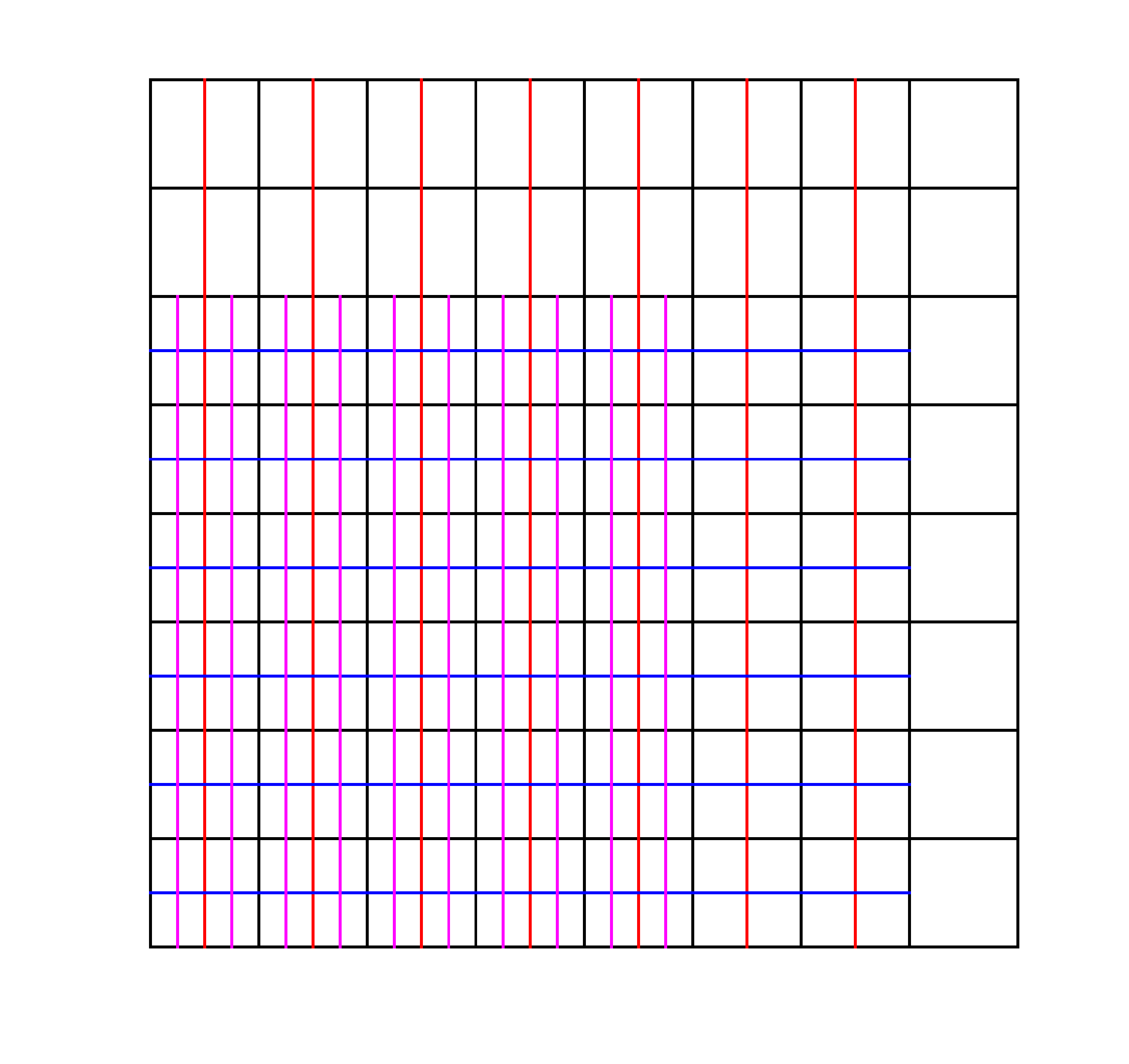}
}
\end{subfigure}
\caption{The fourth generation T-meshes in the first numerical test, and the corresponding B\'ezier meshes, for $p=2$ and $p=3$. For interpretation of the references to color in this figure legend, the reader is referred to the web version of this article} \label{fig3:meshes}
\end{figure}

In the first test the domain is the unit square $\Omega_0 = (0,1)^2$, which is refined near the origin with the following procedure: the level zero mesh is a Cartesian grid, and then for each level we refine, in the T-mesh, the maximal square subregion such that the T-mesh remains \mbox{${\bf p}$-admissible} in the sense of \cite{Morgenstern_Peterseim}. The number of elements of the initial Cartesian grid has been chosen to obtain a similar number of degrees of freedom for each degree. In particular, it is equal to $7^2$, $8^2$ and $10^2$ for biquadratic, bicubic and biquartic splines, respectively. We show in Figure~\ref{fig3:meshes} the refined meshes for the biquadratic and the bicubic case, after all the bisections of the fourth generation.

For the computations with the BPX preconditioner, the inclusion operator $I_k$ can be computed through knot insertion, while the restriction operator is simply its transpose. We have compared the results with two different smoothers $R_k$: Jacobi and symmetric Gauss-Seidel, with one single iteration each. The condition number of the preconditioned system is estimated computing the minimum and maximum eigenvalues with Lanczos' method. The stopping criterion is set to an initial tolerance of $10^{-6}$ for coarse meshes. Table~\ref{tab:square} summarizes the obtained results, where the number of degrees of freedom corresponds to the matrix size, that is, after applying Dirichlet boundary conditions.

As predicted by the theory, the numerical results show that the condition number remains bounded by a constant during $h$-refinement. As already observed in \cite{BHKS13}, the constant deteriorates with the degree $p$, and as explained in \cite{Hofreither2016,Hofreither-Takacs} (see also \cite{Serra-MG-Galerkin-2016}) this is due to the bad behavior of Jacobi and symmetric Gauss-Seidel smoothers with splines. {\Bd The smoother recently proposed in \cite{Hofreither-Takacs} for tensor-product B-splines can be probably adapted to AS T-splines by bisection, although we expect a considerable loss of efficiency of the preconditioner, since the construction of the smoother in \cite{Hofreither-Takacs} strongly relies on the Kronecker tensor-product structure of the matrices.} We also remark that, although the condition number is lower using symmetric Gauss-Seidel as a smoother, the iteration cost in our implementation is higher, thus in our tests the Jacobi smoother performs better than the Gauss-Seidel one.


\begin{table}[ht]
\hspace{-1cm}
\begin{tabular}{|c||c|c|c|c||c|c|c|c||c|c|c|c|}
\hline
 & \multicolumn{4}{|c||}{Biquadratic} & \multicolumn{4}{|c||}{Bicubic} & \multicolumn{4}{|c|}{Biquartic} \\
\hline Levels & Dofs & N.P. & Jac. & GS & Dofs & N.P. & Jac. & GS & Dofs & N.P. & Jac. & GS \\ \hline
2& 85&	14.8&	8.0&	2.7&	137&	48.3&	33.8&	6.6&	234&	373.9&	167.7&	23.9 \\
3& 135&	23.7&	10.6&	3.5&	215&	33.6&	56.7&	7.6&	339&	272.1&	294.5&	35.4 \\
4& 216&	46.5&	14.7&	5.3&	325&	64.6&	74.2&	13.2&	495&	373.1&	393.7&	46.5 \\
5& 344&	74.2&	16.7&	6.0&	496&	75.2&	88.4&	12.9&	705&	299.6&	494.7&	58.8 \\
6& 569&	139.3&	20.6&	7.5&	768&	181.0&	103.6&	19.7&	1047&	378.4&	572.2&	67.0 \\
7& 961&	234.5&	20.4&	7.6&	1233&	217.6&	109.6&	20.4&	1575&	313.8&	642.8&	82.6 \\
8& 1690&	447.7&	23.2&	8.7&	2045&	540.6&	123.6&	24.9&	2505&	737.3&	694.3&	86.0 \\
9& 3042&	798.9&	22.1&	8.5&	3530&	686.1&	124.3&	24.8&	4101&	918.9&	743.1&	103.6 \\
10& 5643&	1560.9&	24.6&	9.4&	6286&	1785.7&	133.2&	27.9&	7071&	2286.9&	770.1&	108.7 \\
11& 10643&	2907.2&	23.1&	9.1&	11539&	2371.1&	133.9&	27.4&	12531&	3030.4&	805.1&	118.1 \\
12& 20444&	5770.1&	25.5&	9.8&	21639&	6382.3&	140.8&	29.9&	23037&	7868.2&	816.5&	121.1 \\
13& 39652&	11040.6&	24.0&	9.4&	41340&	8729.9&	140.4&	29.4&	43137&	10855.7&	843.2&	127.1 \\
14& 77677&	22109.3&	26.2&	9.9&	79952&	23983.8&	145.9&	31.2&	82539&	28931.9&	846.6&	 128.1 \\
15& 152949&	42979.9&	24.6&	9.6&	156197&	33390.6&	145.0&	30.6&	159567&	40886.2&	867.2&	 132.7 \\ \hline
\end{tabular}
\caption{Condition numbers for the square domain. N.P: no preconditioner. Jac: Jacobi smoother. GS: Gauss-Seidel smoother} \label{tab:square}
\end{table}

\subsubsection{Curved L-shaped domain}
\begin{figure}[ht]
\centerline{\includegraphics[width=0.5\textwidth,trim=3cm 2cm 2cm 1cm, clip]{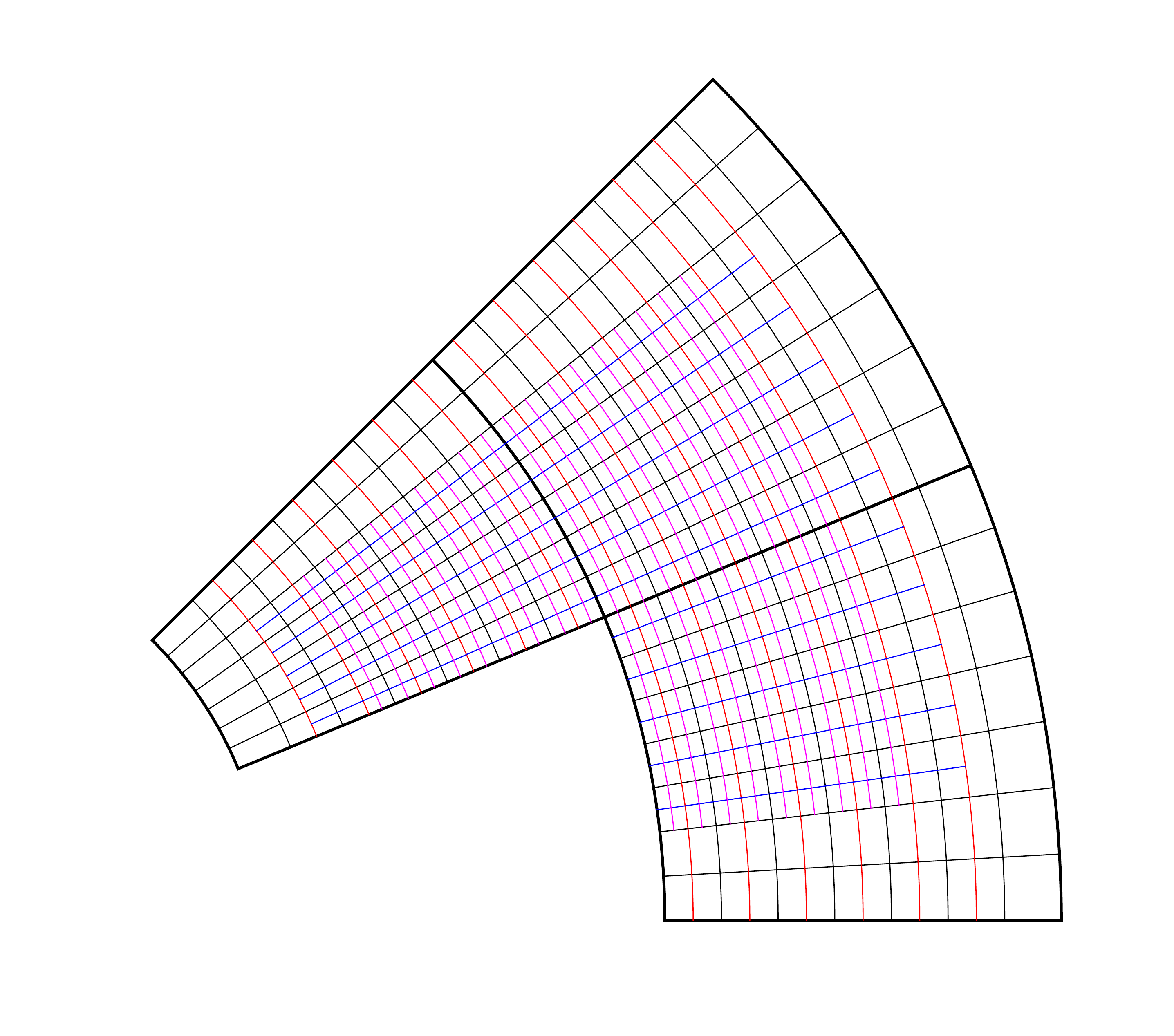}}
\caption{B\'ezier mesh of the curved L-shaped domain with four levels, biquadratic case. For interpretation of the references to color in this figure legend, the reader is referred to the web version of this article} \label{fig:curvedL}
\end{figure}

As a second test we have considered the curved L-shaped domain in Figure~\ref{fig:curvedL}. The domain is defined with three patches, each one the image of the unit square through a different parametrization ${\bf F}$, in such a way that the reentrant corner is always the image of the origin. We have constructed, in the parametric domain, the same meshes as in the previous test, which are then mapped to the physical domain through ${\bf F}$. The relative orientation of the three patches guarantees that the meshes match on the interfaces, and basis functions on the interfaces are glued together with $C^0$ continuity.

We have run the same kind of tests as before, the results are summarized in Table~\ref{tab:curvedL}. The numerical results show that the presence of the parametrization does not greatly affect the condition number, as long as the parametrization remains regular. These results are in agreement with those in \cite{BHKS13} for tensor-product B-splines.

\begin{table}
\hspace{-1cm}
\begin{tabular}{|c||c|c|c|c||c|c|c|c||c|c|c|c|}
\hline
 & \multicolumn{4}{|c||}{Biquadratic} & \multicolumn{4}{|c||}{Bicubic} & \multicolumn{4}{|c|}{Biquartic} \\
\hline Levels & Dofs & N.P. & Jac. & GS & Dofs & N.P. & Jac. & GS & Dofs & N.P. & Jac. & GS \\ \hline
2&	275&	60.7&	10.2&	2.8&	436&	122.6&	47.1&	6.7&	735&	1060.8&	233.4&	23.9 \\
3&	430&	80.0&	16.2&	4.0&	676&	147.9&	82.0&	9.9&	1057&	744.6&	394.5&	36.5 \\
4&	682&	202.4&	22.6&	5.6&	1016&	368.5&	113.8&	13.4&	1537&	1048.1&	581.4&	48.9 \\
5&	1074&	270.5&	27.5&	6.6&	1538&	478.7&	145.9&	16.7&	2177&	870.5&	756.6&	61.6 \\
6&	1764&	654.7&	32.8&	8.0&	2370&	1117.5&	172.9&	20.1&	3221&	1963.0&	910.3&	74.5 \\
7&	2954&	914.5&	36.3&	8.5&	3780&	1509.3&	198.4&	22.7&	4821&	2506.6&	1047.1&	86.0 \\
8&	5168&	2203.6&	41.1&	9.6&	6244&	3510.3&	220.3&	25.6&	7641&	5633.6&	1168.9&	97.9 \\
9&	9250&	3238.5&	43.6&	9.7&	10726&	5013.1&	240.2&	27.4&	12457&	7710.7&	1276.9&	106.5 \\
10&	17104&	7864.6&	47.8&	10.6&	19046&	11862.7&	258.5&	29.4&	21421&	17773.4&	1375.3&	116.2 \\
11&	32154&	12027.0&	49.8&	10.5&	34856&	17793.4&	275.0&	30.4&	37853&	25944.3&	1461.5&	 122.0 \\
12&	61656&	29424.9&	53.5&	11.3&	65256&	42832.1&	290.1&	31.9&	69473&	61381.7&	1541.2&	 129.2 \\
13&	119378&	46178.5&	54.9&	11.1&	124458&	66460.2&	303.5&	32.4&	129873&	93799.5&	1610.5&	 132.5 \\
\hline
\end{tabular}
\caption{Condition numbers for the curved L-shaped domain. N.P: no preconditioner. Jac: Jacobi smoother. GS: Gauss-Seidel smoother} \label{tab:curvedL}
\end{table}

\subsection{Comparison with other variants}
\subsubsection{Micro test}
Our implementation is based on tensor-product spaces, therefore it is not general enough to test the preconditioner using the \emph{micro} technique, that is, defining the decomposition spaces ${\cal V}_k$ as in \eqref{subspace_Vk} adding one edge at each step. However, in order to better understand how the \emph{micro} decomposition could work we have tried a modified version of the method: instead of adding all the bisections of the same generation at the same time, as in the \emph{macro} technique, we group aligned bisection edges of the same generation into one single level. For instance, in the example of Figure~\ref{fig:Tmesh-square-p3} for the bicubic case, and recalling that we start from generation zero, the first, second and third generations contain seven, six and ten new levels, respectively.

{\Bd In Table~\ref{tab:micro} we show the condition numbers of the system obtained after all the lines of the same generation have been obtained, that can be compared with the ones obtained with the macro technique in Table~\ref{tab:square}. Notice that we do not show the number of degrees of freedom and the condition number of the unpreconditioned system, since they are the same as in Table~\ref{tab:square}. Although the condition number remains bounded, as predicted by the theory, the results show a loss of efficiency of the preconditioner. This seems to be a consequence of the behavior of the constants in the inequalities \eqref{nonoverlapping_consequence}, that depend on the overlaps of the sets $\omega_k$ (and $\tilde \omega_k$) of the same generation, see also Remark~\ref{rem:constants}. Gathering all the bisections of the same generation together, as in the macro technique, reduces the number of overlaps and the value of these constants, while considering the bisections separately, as in the micro technique, increases the number of overlaps. We expect even a worse behavior for the ``pure'' micro decomposition of Section~\ref{sec:micro}. We also remark that this behavior does not appear in the case of finite elements \cite{CNX}, because of the reduced support of $C^0$ finite element functions.
}


\begin{table}
\centering
\begin{tabular}{|c|c|c||c|c|c||c|c|c|}
\hline
\multicolumn{3}{|c||}{\Bd Biquadratic} & \multicolumn{3}{|c||}{Bicubic} & \multicolumn{3}{|c|}{\Bd Biquartic} \\
\hline
Levels & Jac & GS & Levels & Jac & GS & Levels & Jac & GS \\
\hline
6 (2)  & 13.0 & 4.6 & 8 (2)& 57.8&	13.2 & 9 (2) & 346.3 & 53.6 \\
11 (3)  & 22.3 & 8.6 & 14 (3)&	112.6&	23.7 & 16 (3) & 653.6 & 91.2 \\
20 (4) & 28.7 & 10.0 & 24 (4)&	145.9&	31.1 & 28 (4) & 877.8 & 138.3 \\
28 (5)  & 34.3 & 14.8 & 33 (5)&	182.5&	42.9 & 38 (5) & 1097.8 & 166.9 \\
43 (6) & 41.4 & 14.8 & 49 (6)&	198.5&	45.9 & 56 (6) & 1220.7 & 214.4\\
57 (7) & 44.0 & 19.1 & 64 (7)&	220.6&	56.2 & 72 (7) & 1379.3 & 234.5 \\
84 (8) & 48.4 & 18.6 & 92 (8)&	232.6&	56.6 & 102 (8) & 1504.3 & 266.3 \\
110 (9) & 47.9 & 22.3 & 119 (9)& 247.2&	67.1 & 130 (9) & 1604.3 & 285.0 \\
161 (10) & 51.8 & 20.6 & 171 (10) & 252.2 & 63.5 & 184 (10) & 1704.0 & 303.2\\
211 (11) & 50.3 & 23.6 & 222 (11) & 260.7 & 72.2 & 236 (11) & 1758.4 & 316.8 \\
310 (12) & 54.0 & 21.4 & 322 (12) & 264.1 & 66.3 & 338 (12) & 1826.7 & 326.1 \\
\hline
\end{tabular}
\caption{Condition number obtained with the \emph{micro} technique, after adding all the lines of the same generation (between parentheses)} \label{tab:micro}
\end{table}



\subsubsection{An alternative refinement}
\begin{figure}[ht]
\includegraphics[width=0.45\textwidth]{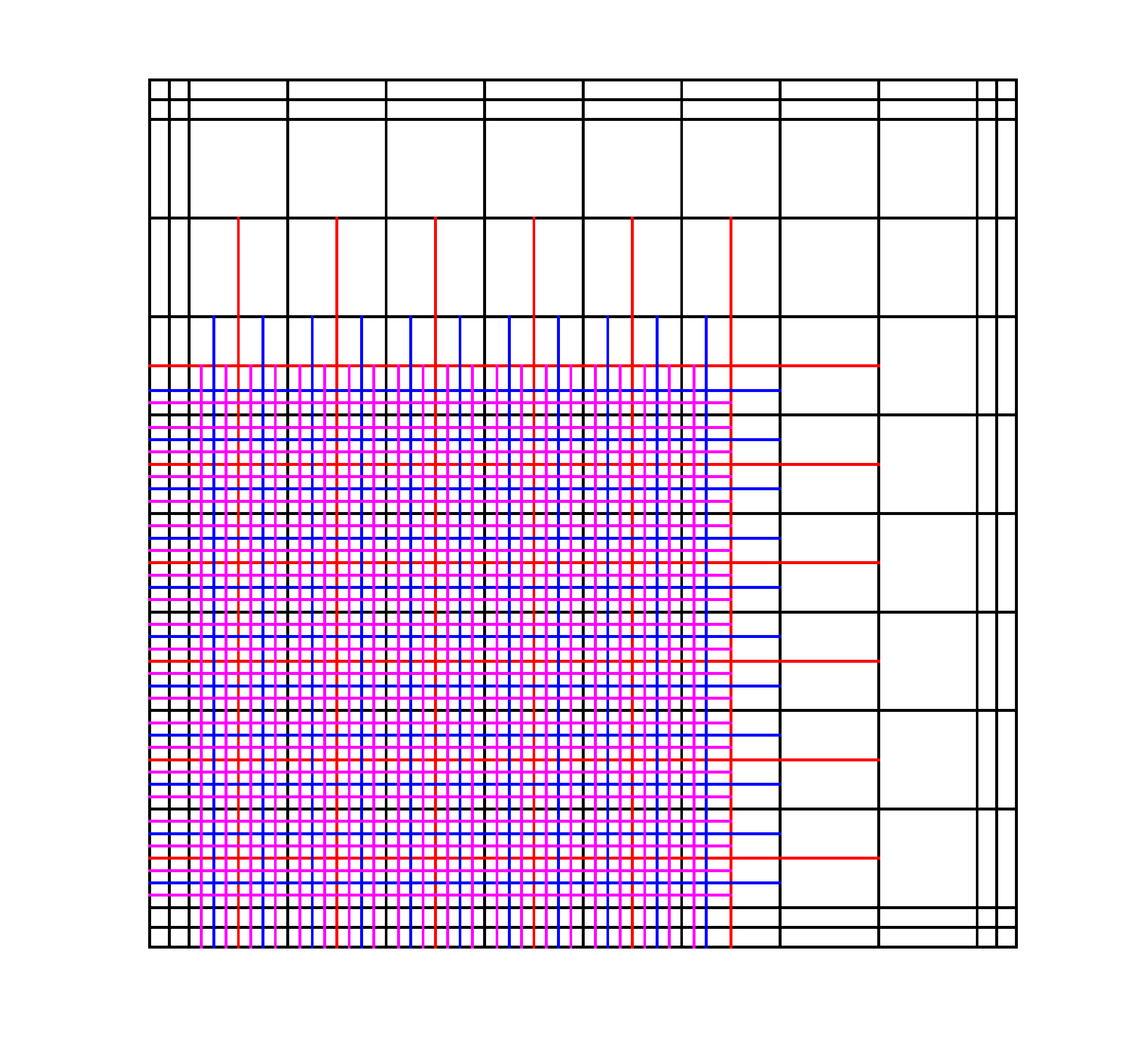}
\includegraphics[width=0.45\textwidth]{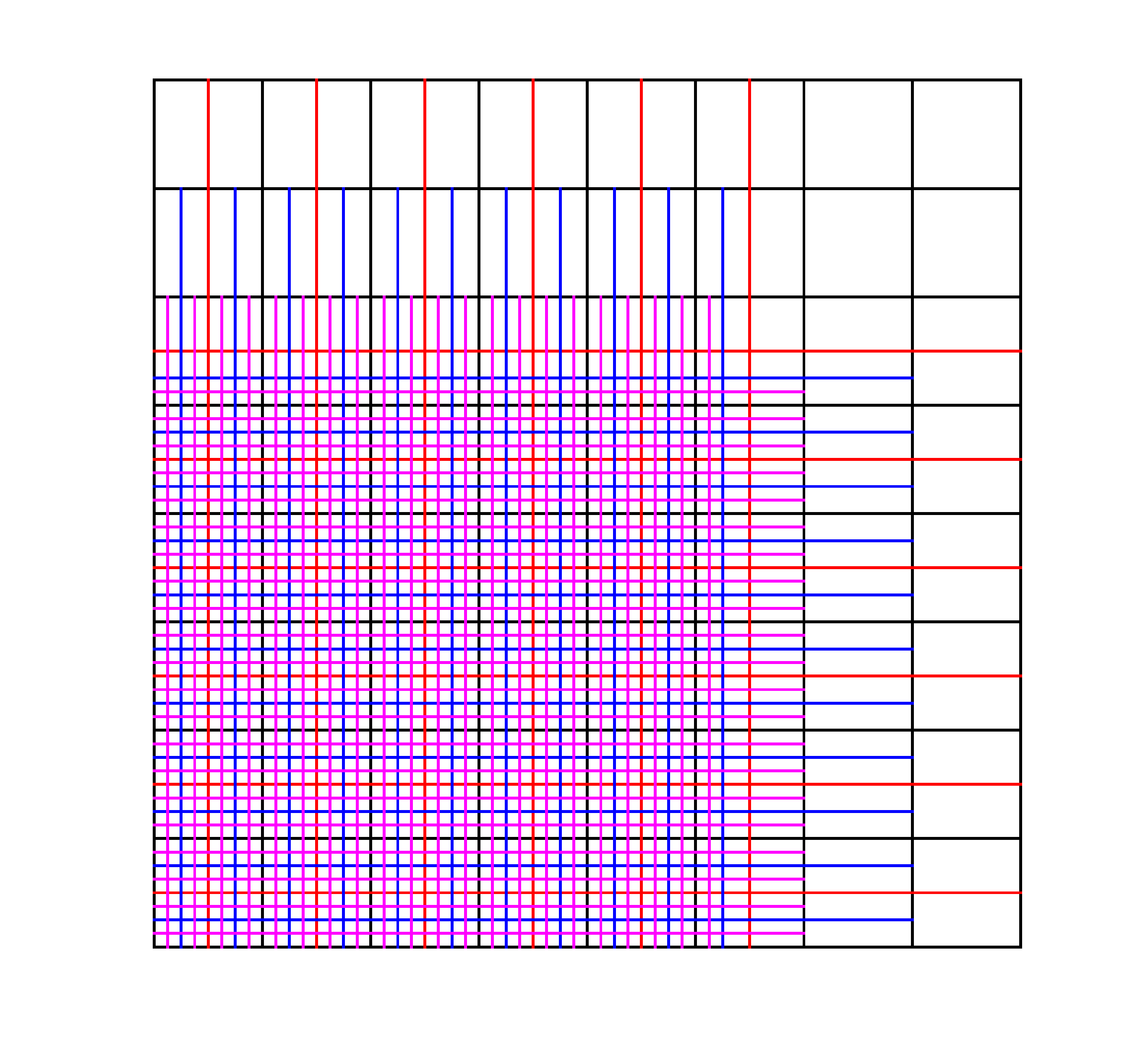}
\caption{The alternative refinement for the biquadratic case: T-mesh (left) and corresponding B\'ezier mesh (right). For interpretation of the references to color in this figure legend, the reader is referred to the web version of this article} \label{fig5:meshes}
\end{figure}

In all the previous tests the refinement is performed as in \cite{Morgenstern_Peterseim} to obtain ${\bf p}$-admissible meshes, that is, the elements are refined alternatively in the vertical and the horizontal direction, depending on their generation. {\Bd For this new test we have tried a different kind of refinement, suitable to refine towards a corner, that preserves the multilevel structure and the analysis-suitability of the mesh. This refinement gives a higher proportion of square elements with respect to ${\bf p}$-admissible meshes, although more elongated elements can also appear.} Starting from a Cartesian grid, the elements of the T-mesh in a square subregion next to the corner are split in four, bisecting them vertically and horizontally simultaneously. Then, some other elements have to be bisected, either horizontally or vertically, to maintain the analysis-suitability property. To simplify the computations, in this second step we refine together a set of aligned elements, although this is not really necessary to maintain the analysis-suitable property. All the functions added or modified during the refinement of the square region are considered to have the same generation, independently of whether the refinement is horizontal, vertical or both. A detail of the meshes obtained after three refinement steps, for the biquadratic case, is shown in Figure~\ref{fig5:meshes}. Lines added at the same step are drawn with the same color.

For the numerical tests we use the \emph{macro} technique as in Section~\ref{sec:macro-tests}, that is, the space ${\cal W}_k$ is defined after adding all the lines of the same generation. In Table~\ref{tab:old-refinement} we show the results obtained applying the BPX preconditioner with this kind of refinement. The results show a better behavior of the BPX preconditioner compared to the meshes in Section~\ref{sec:macro-tests}, and in fact the numbers are very similar to the ones obtained for tensor-product B-splines in \cite{BHKS13}. Indeed, this seems a good motivation to look for alternative refinement techniques for T-splines, that are less restrictive than the one in \cite{Morgenstern_Peterseim}. This may be possible after the definition of truncated T-splines \cite{TT-splines}, for which has been proved linear independence with a less restrictive condition than the one of analysis-suitable T-splines.

\begin{table}
\hspace{-1cm}
\begin{tabular}{|c||c|c|c|c||c|c|c|c||c|c|c|c|}
\hline
 & \multicolumn{4}{|c||}{Biquadratic} & \multicolumn{4}{|c||}{Bicubic} & \multicolumn{4}{|c|}{Biquartic} \\
\hline Levels & Dofs & N.P. & Jac. & GS & Dofs & N.P. & Jac. & GS & Dofs & N.P. & Jac. & GS \\ \hline
2& 184&	23.4&	6.4&	2.4&	213&	33.5&	30.0&	5.9&	244&	243.3&	134.6&	19.8 \\
3& 569&	103.5&	10.2&	3.6&	620&	126.6&	48.8&	9.0&	673&	390.8&	285.8&	34.8 \\
4& 1934&	402.8&	11.3&	4.4&	2027&	608.9&	60.5&	12.0&	2122&	853.1&	368.3&	47.0 \\
5& 7059& 1530.0& 12.0&	4.9&	7234&	2542.9&	65.2&	14.0&	7411&	4318.7&	409.0&	57.2 \\
6& 26904& 5999.2& 12.3&	5.0&	27241&	10270.6& 67.3&	15.0&	27580&	18742.2& 423.4&	63.5 \\
7& 104989& 23738.0& 12.5& 5.1&	105648&	41233.4& 68.9&	15.6&	106309&	77464.3& 427.3&	67.0 \\
\hline
\end{tabular}
\caption{Condition numbers with the alternative refinement. N.P: no preconditioner. Jac: Jacobi smoother. GS: Gauss-Seidel smoother} \label{tab:old-refinement}
\end{table}

\section{Conclusions}

We have presented the analysis of an additive multilevel preconditioner for T-splines. The T-meshes are assumed to be ${\bf p}$-admissible, as in the construction by \cite{Morgenstern_Peterseim}, which guarantees analysis-suitability of the T-splines. The optimality of the preconditioner under $h$-refinement is proved following the ideas by \cite{CNX}, and relies on a stable space decomposition (Theorem~\ref{stabledecomp_ASTmeshes}) and the strengthened Cauchy-Schwarz inequality (Theorem~\ref{SCS_AST}). In the method by \cite{CNX} a new level is defined each time we add one edge to the mesh, and only the functions modified by this edge are added to the level. We have studied an alternative construction, where the new level is defined after adding all the edges of the same generation, and considering the newly added and modified functions as having the same level. The theoretical analysis for this alternative decomposition follows from the previous one, but the numerical results show a better performance. \Bd We also remark that, since the refinement strategy in \cite{Morgenstern_Peterseim} has been generalized in \cite{Morgenstern} to the three dimensional case with trivariate odd degree AS T-tsplines, our proposed preconditioners can be also generalized to the three dimensional case accordingly.
\Bk


\section*{Acknowledgments}
The authors would like to thank Prof. Annalisa Buffa for many fruitful discussions on the subject of the paper. The work of Durkbin Cho was supported by Basic Science Research Program through the National Research Foundation of Korea (NRF) funded by the Ministry of Science, ICT \& Future Planning (2015R1A1A1A05001109). The work of Rafael V\'azquez was partially supported by the European Research Council through the FP7 ERC Consolidator Grant no.~616563 HIGEOM (PI: Giancarlo Sangalli) and the FP7 ERC Advanced Grant no.~694515 CHANGE (PI: Annalisa Buffa), and by European Union's Horizon 2020 research and innovation programme through the grant no.~680448 CAxMan. This support is gratefully acknowledged.

This is a pre-copyedited, author-produced PDF of an article accepted for publication in IMA Journal of Numerical Analysis following peer review. The version of record is available online at: \url{https://doi.org/10.1093/imanum/dry032}.

\appendix
\section{Mathematical proofs of auxiliary results}
\renewcommand{\thesection}{A}

\subsection{Difference of T-spline projectors applied to coarse B-splines}\label{appendixA1}
The first auxiliary result was used in the proof of Theorem~\ref{stabledecomp_ASTmeshes}. Letting $\bar{v}_\ell :=({\bf \Pi}_{\bf p, \Xi^\ell}-{\bf \Pi}_{\bp, {\bf \Xi}^{\ell-1}})v \in S_\bp({\bf \Xi}^\ell),\ 0\le \ell \le L$, we need to prove that for $l\le \ell_k-1$
\[
({\bf \Pi}_{\bf p}^{{\cal T}_k}-{\bf \Pi}_{\bf p}^{{\cal T}_{k-1}})\bar{v}_l=0,
\]
which is a consequence of the following lemma.
\begin{lemm}
Let ${\cal T}_k = {\cal T}_{k-1} + b_{\tau_{k-1}}$ be a $\bp$-admissible T-mesh generated by bisection of $\tau_{k-1}$, and let $\ell_k = g(\tau_{k-1})+1$. Then it holds that
\[
({\bf \Pi}_{\bf p}^{{\cal T}_k}-{\bf \Pi}_{\bf p}^{{\cal T}_{k-1}}) v = 0
\]
for any function $v \in S_\bp({\bf \Xi}^\ell)$, with $\ell < \ell_k$.
\end{lemm}
{\it Proof.}
Since $S_\bp({\bf \Xi}^i) \subset S_\bp ({\bf \Xi}^j)$ for $i \le j$, we only need to prove the result for a function $v \in S_\bp({\bf \Xi}^{\ell_k - 1})$. Let us first define, similarly to $\Phi_k$, the set
\[
\Psi_{k-1} := \{B_{\bA,\bp} : \bA \in {\cal A}_\bp({\cal T}_{k-1}) \} \setminus \{B_{\bA,\bp} : \bA \in {\cal A}_\bp({\cal T}_{k}) \},
\]
which is the collection of functions removed after the bisection $b_{\tau_{k-1}}$. Since the bisection only affects basis functions in $\Psi_{k-1}$ (removed) and $\Phi_k$ (added), from the definition of the dual functionals it is clear that
\[
({\bf \Pi}_{\bf p}^{{\cal T}_k}-{\bf \Pi}_{\bf p}^{{\cal T}_{k-1}}) v = \sum_{\bA : B_{\bA,\bp} \in \Phi_k} (\lambda_{\bA,\bp} v) B_{\bA,\bp} - \sum_{\bA : B_{\bA,\bp} \in \Psi_{k-1}} (\lambda_{\bA,\bp} v) B_{\bA,\bp}.
\]

From the nestedness of the T-spline spaces $S({\cal A}_\bp({\cal T}_{k-1})) \subset S({\cal A}_\bp({\cal T}_k))$ (see \cite[Corollary~5.8]{Morgenstern_Peterseim}), it holds that $\mbox{span} (\Psi_{k-1}) \subset \mbox{span} (\Phi_{k})$, and by the local linear independence of analysis-suitable T-splines, we have
\[
\bigcup_{B_{\bA,\bp} \in \Psi_{k-1}} \mbox{supp} (B_{\bA,\bp}) \subset \bigcup_{B_{\bA,\bp} \in \Phi_{k}} \mbox{supp} (B_{\bA,\bp}) = \omega_k.
\]
Moreover, from Lemma~\ref{lemma:lqu}, and in particular from \eqref{hsk-U}-\eqref{vsk-U}, we know that functions in $\Phi_k$ see all the lines of the uniform Cartesian mesh ${\cal T}^u_{\ell_k-1}$ in the vicinity of $\tau_{k-1}$, and using the previous property for the supports, we know that the same is true for functions in $\Psi_{k-1}$.

Let us define, starting from the uniform Cartesian mesh ${\cal T}'_0 = {\cal T}^u_{\ell_k-1}$ the auxiliary family of $\bp$-admissible T-meshes
\[
{\cal T}'_{j+1} = {\cal T}'_{j} + b_{\tau_{j + k_0}}, \quad j = 0, \ldots, k-k_0-1,
\]
where $k_0$ is the minimum integer such that $g(\tau_{k_0}) = g(\tau_{k-1})$ in the construction \eqref{eq:T-by-bisection} of the T-mesh ${\cal T}_{k+1}$. In other words,  starting from the Cartesian grid of level $\ell_k-1$, we perform the bisections of the same generation until we reach $b_{\tau_{k-1}}$. From what we have just seen, it is clear that $\Phi'_{k-k_0} = \Phi_k$ and $\Psi'_{k-k_0-1} = \Psi_{k-1}$, and by the nestedness of $\bp$-admissible T-splines it holds $S_\bp({\bf \Xi}^{\ell_k-1}) = S({\cal A}_\bp({\cal T}'_0)) \subset S({\cal A}_\bp({\cal T}'_{j}))$ for any $j = 0, \ldots, k-k_0$. Thus we have, for $v \in S_\bp({\bf \Xi}^{\ell_k-1})$
\[
0 = \Big( {\bf \Pi}_\bp^{{\cal T}'_{k-k_0}} - {\bf \Pi}_\bp^{{\cal T}'_{k-k_0-1}} \Big) v = \sum_{\bA : B_{\bA,\bp} \in \Phi'_{k-k_0}} (\lambda_{\bA,\bp} v) B_{\bA,\bp} - \sum_{\bA : B_{\bA,\bp} \in \Psi'_{k-k_0-1}} (\lambda_{\bA,\bp} v) B_{\bA,\bp},
\]
which ends the proof. \hfill $\square$

\subsection{Appropriate discrete norms}\label{appendixA2}
The second auxiliary result, that we present in this subsection, is $\|u_{\rm B}\|_{L^2(\Gamma)} \lesssim  \|u_j\|_{L^2(Q^i)}$ in the proof of Lemma~5.3.\\

Let $Q^i \in { \cal M}^{\cal B}_i$ be a B\'ezier element of the Cartesian grid of generation $i$. For $j \ge i$ we define the set of indices
\[
{\bf I}^j(Q^i):= \{ {\bf j} \in {\bf I}^j : {\rm supp} (B_{\bf j, p}) \cap Q^i \ne \emptyset \},
\]
which is also the union of indices ${\bf I}^j_{\rm in}$ and ${\bf I}^j_{\rm B}$ defined as in the proof of Lemma \ref{SCS_quasiuniform}, namely, ${\bf I}^j(Q^i)= {\bf I}^j_{\rm in} \cup {\bf I}^j_{\rm B}$. By local linear independence of B-splines, a function $z \in S_\bp({\bf \Xi}^j)$ restricted to $Q^i$ can be written as
\[
z|_{Q^i} = \sum_{{\bf j} \in {\bf I}^j(Q^i)} c_{\bf j} B_{\bf j, p}.
\]
Let also $Q^j \in {\cal M}^{\cal B}_j$ such that $Q^j \subset Q^i$, and define ${\bf I}^j(Q^j) := \{ {\bf j} \in {\bf I}^j : Q^j \subset {\rm supp}(B_{\bf j, p}) \}$. Obviously,
${\bf I}^j(Q^j) \subset {\bf I}^j(Q^i)$.
We define a local discrete norm $|\cdot|_{Q^j}$ as follows: for any $z\in  S_\bp({\bf \Xi}^j)|_{Q^i}$,
\[
|z|^2_{Q^j}:=\left( \max_{{\bf j} \in {\bf I}^j(Q^j)} |c_{\bf j}|^2\right) h^2_j,
\]
and also define the global discrete norm $|\cdot|_{Q^i}$ on $Q^i$ as follows: for any $z\in  S_\bp({\bf \Xi}^j)|_{Q^i}$,
\[
|z|_{Q^i}^2 := \sum_{{\bf j}\in {\bf I}^j(Q^i)} |c_{\bf j}|^2 h_j^2.
\]

\begin{proposition}\label{discrete_norm} Let $Q^i\in {\cal M}^{\cal B}_i$ and let, for $j \ge i$, $Q^j\in {\cal M}^{\cal B}_j$ such that $Q^j \subset Q^i$. Then for any $z\in  S_\bp({\bf \Xi}^j)|_{Q^i}$ it holds
\[
\|z\|_{L^2(Q^j)}^2 \approx |z|^2_{Q^j},
\]
and
\[
\|z\|_{L^2(Q^i)}^2 \approx  |z|^2_{Q^i}.
\]
\end{proposition}
{\it Proof.} The proof can be found in \cite[Proposition~5.1 and Corollary~5.1]{{BePa_BDDC}}. \hfill $\square$

By Proposition~\ref{discrete_norm}, it is clear that
\[
\|u_{\rm B}\|_{L^2(\Gamma)} \lesssim \|u_j\|_{L^2(Q^i)}.
\]

\bibliographystyle{siam}
\bibliography{biblio}

\def\cprime{$'$}
\begin{thebibliography}{10}

\bibitem{Bazilevs_Beirao_Cottrell_Hughes_Sangalli}
{\sc Y.~Bazilevs, L.~Beir{\~a}o~da Veiga, J.~A. Cottrell, T.~J.~R. Hughes, and
  G.~Sangalli}, {\em Isogeometric analysis: approximation, stability and error
  estimates for {$h$}-refined meshes}, Math. Models Methods Appl. Sci., 16
  (2006), pp.~1031--1090.

\bibitem{Bazilervs_Calo_Cottrell_Evans}
{\sc Y.~Bazilevs, V.~Calo, J.~A. Cottrell, J.~A. Evans, T.~J.~R. Hughes,
  S.~Lipton, M.~Scott, and T.~Sederberg}, {\em Isogeometric analysis using
  {T}-splines}, Comput. Methods Appl. Mech. Engrg., 199 (2010), pp.~229 -- 263.

\bibitem{BPSWZ}
{\sc L.~Beir\~ao~da Veiga, L.~F. Pavarino, S.~Scacchi, O.~B. Widlund, and
  S.~Zampini}, {\em Isogeometric {BDDC} preconditioners with deluxe scaling},
  SIAM J. Sci. Comput., 36 (2014), pp.~A1118--A1139.

\bibitem{BBCS12}
{\sc L.~Beir{\~a}o~da Veiga, A.~Buffa, D.~Cho, and G.~Sangalli}, {\em
  Analysis-{S}uitable {T}-splines are {D}ual-{C}ompatible}, Comput. Methods
  Appl. Mech. Engrg., 249--252 (2012), pp.~42 -- 51.

\bibitem{BBSV2013}
{\sc L.~Beir{\~a}o~da Veiga, A.~Buffa, G.~Sangalli, and R.~V\'azquez}, {\em
  Analysis-suitable {T}-splines of arbitrary degree: Definition, linear
  independence and approximation properties}, Math. Models Methods Appl. Sci.,
  23 (2013), pp.~1979--2003.

\bibitem{IGA-acta}
{\sc L.~Beir{\~a}o~da Veiga, A.~Buffa, G.~Sangalli, and R.~V{\'a}zquez}, {\em
  Mathematical analysis of variational isogeometric methods}, Acta Numer., 23
  (2014), pp.~157--287.

\bibitem{BePa_BDDC}
{\sc L.~Beir{\~a}o~da Veiga, D.~Cho, L.~F. Pavarino, and S.~Scacchi}, {\em
  {BDDC} preconditioners for isogeometric analysis}, Math. Models Methods Appl.
  Sci., 23 (2013), pp.~1099 -- 1142.

\bibitem{BePa13}
\leavevmode\vrule height 2pt depth -1.6pt width 23pt, {\em Isogeometric
  {S}chwarz preconditioners for linear elasticity systems}, Comput. Methods
  Appl. Mech. Engrg., 253 (2013), pp.~439 -- 454.

\bibitem{BPWX}
{\sc J.~H. Bramble, J.~E. Pasciak, J.~P. Wang, and J.~Xu}, {\em Convergence
  estimates for multigrid algorithms without regularity assumptions}, Math.
  Comp., 57 (1991), pp.~23--45.

\bibitem{BPX1}
{\sc J.~H. Bramble, J.~E. Pasciak, and J.~Xu}, {\em Parallel multilevel
  preconditioners}, Math. Comp., 55 (1990), pp.~1--22.

\bibitem{Brenner_Scott_1994}
{\sc S.~C. Brenner and L.~R. Scott}, {\em The mathematical theory of finite
  element methods}, vol.~15 of Texts in Applied Mathematics, Springer-Verlag,
  New York, 1994.

\bibitem{Bressan_Buffa_Sangalli}
{\sc A.~Bressan, A.~Buffa, and G.~Sangalli}, {\em Characterization of
  analysis-suitable {T}-splines}, Comput. Aided Geom. Design, 39 (2015), pp.~17
  -- 49.

\bibitem{BHKS13}
{\sc A.~Buffa, H.~Harbrecht, A.~Kunoth, and G.~Sangalli}, {\em
  {BPX}-preconditioning for isogeometric analysis}, Comput. Methods Appl. Mech.
  Engrg., 265 (2013), pp.~63 -- 70.

\bibitem{CNX}
{\sc L.~Chen, R.~H. Nochetto, and J.~Xu}, {\em Optimal multilevel methods for
  graded bisection grids}, Numer. Math., 120 (2012), pp.~1--34.

\bibitem{IGA-book}
{\sc J.~A. Cottrell, T.~J.~R. Hughes, and Y.~Bazilevs}, {\em Isogeometric
  {A}nalysis: toward integration of {CAD} and {FEA}}, John Wiley \& Sons, 2009.

\bibitem{DeBoor}
{\sc C.~de~Boor}, {\em A practical guide to splines}, vol.~27 of Applied
  Mathematical Sciences, Springer-Verlag, New York, revised~ed., 2001.

\bibitem{dFRV11}
{\sc C.~de~Falco, A.~Reali, and R.~V\'azquez}, {\em Geo{PDE}s: a research tool
  for {I}sogeometric {A}nalysis of {PDE}s}, Adv. Engrg. Softw., 42 (2011),
  pp.~1020--1034.

\bibitem{LR-splines}
{\sc T.~Dokken, T.~Lyche, and K.~F. Pettersen}, {\em Polynomial splines over
  locally refined box-partitions}, Comput. Aided Geom. Design, 30 (2013),
  pp.~331--356.

\bibitem{Serra-MG-Galerkin-2016}
{\sc M.~Donatelli, C.~Garoni, C.~Manni, S.~Serra-Capizzano, and H.~Speleers},
  {\em Robust and optimal multi-iterative techniques for {I}g{A} {G}alerkin
  linear systems}, Comput. Methods Appl. Mech. Engrg., 284 (2015),
  pp.~230--264.

\bibitem{Manni_MG}
\leavevmode\vrule height 2pt depth -1.6pt width 23pt, {\em Symbol-{B}ased
  {M}ultigrid {M}ethods for {G}alerkin {B}-{S}pline {I}sogeometric {A}nalysis},
  SIAM J. Numer. Anal., 55 (2017), pp.~31--62.

\bibitem{Dorfel_Juttler_Simeon}
{\sc M.~D\"{o}rfel, B.~J\"{u}ttler, and B.~Simeon}, {\em Adaptive isogeometric
  analysis by local $h$-refinement with {T}-splines}, Comput. Methods Appl.
  Mech. Engrg., 199 (2010), pp.~264 -- 275.

\bibitem{Gahalaut_MG}
{\sc K.~Gahalaut, J.~Kraus, and S.~Tomar}, {\em Multigrid methods for
  isogeometric discretization}, Comput. Methods Appl. Mech. Engrg., 253 (2013),
  pp.~413 -- 425.

\bibitem{Hofreither-Takacs}
{\sc C.~Hofreither and S.~Takacs}, {\em Robust multigrid for isogeometric
  analysis based on stable splittings of spline spaces}, SIAM J. Num. Anal., 55
  (2017), pp.~2004--2024.

\bibitem{Hofreither2016}
{\sc C.~Hofreither, S.~Takacs, and W.~Zulehner}, {\em A robust multigrid method
  for {I}sogeometric {A}nalysis in two dimensions using boundary correction},
  Comput. Methods Appl. Mech. Engrg., 316 (2017), pp.~22 -- 42.

\bibitem{Hofreither2015_1}
{\sc C.~Hofreither and W.~Zulehner}, {\em Mass smoothers in geometric multigrid
  for isogeometric analysis}, in Curves and surfaces, vol.~9213 of Lecture
  Notes in Comput. Sci., Springer, Cham, 2015, pp.~272--279.

\bibitem{Hofreither2015_2}
\leavevmode\vrule height 2pt depth -1.6pt width 23pt, {\em Spectral analysis of
  geometric multigrid methods for isogeometric analysis}, in Numerical methods
  and applications, vol.~8962 of Lecture Notes in Comput. Sci., Springer, Cham,
  2015, pp.~123--129.

\bibitem{Hughes_Cottrell_Bazilevs}
{\sc T.~J.~R. Hughes, J.~A. Cottrell, and Y.~Bazilevs}, {\em Isogeometric
  analysis: {CAD}, finite elements, {NURBS}, exact geometry and mesh
  refinement}, Comput. Methods Appl. Mech. Engrg., 194 (2005), pp.~4135--4195.

\bibitem{Johannessen2013}
{\sc K.~A. Johannessen, T.~Kvamsdal, and T.~Dokken}, {\em Isogeometric analysis
  using {LR} {B}-splines}, Comput. Methods Appl. Mech. Engrg., 269 (2014),
  pp.~471 -- 514.

\bibitem{KPJ_FETIDP}
{\sc S.~K. Kleiss, C.~Pechstein, B.~J\"uttler, and S.~Tomar}, {\em
  I{ETI}---isogeometric tearing and interconnecting}, Comput. Methods Appl.
  Mech. Engrg., 247/248 (2012), pp.~201--215.

\bibitem{Li-Scott}
{\sc X.~Li and M.~A. Scott}, {\em Analysis-suitable {T}-splines:
  {C}haracterization, refineability, and approximation}, Math. Models Methods
  Appl. Sci., 24 (2014), pp.~1141--1164.

\bibitem{LZSHS12}
{\sc X.~Li, J.~Zheng, T.~Sederberg, T.~Hughes, and M.~Scott}, {\em On linear
  independence of {T}-spline blending functions}, Comput. Aided Geom. Design,
  29 (2012), pp.~63 -- 76.

\bibitem{Morgenstern}
{\sc P.~Morgenstern}, {\em Globally structured three-dimensional
  analysis-suitable {T}-splines: definition, linear independence and
  {$m$}-graded local refinement}, SIAM J. Numer. Anal., 54 (2016),
  pp.~2163--2186.

\bibitem{Morgenstern_Peterseim}
{\sc P.~Morgenstern and D.~Peterseim}, {\em Analysis-suitable adaptive {T}-mesh
  refinement with linear complexity}, Comput. Aided Geom. Design, 34 (2015),
  pp.~50--66.

\bibitem{PS_FETIDP}
{\sc L.~F. Pavarino and S.~Scacchi}, {\em Isogeometric block {FETI}-{DP}
  preconditioners for the {S}tokes and mixed linear elasticity systems},
  Comput. Methods Appl. Mech. Engrg., 310 (2016), pp.~694--710.

\bibitem{Sangalli_Tani}
{\sc G.~Sangalli and M.~Tani}, {\em Isogeometric preconditioners based on fast
  solvers for the {S}ylvester equation}, SIAM J. Sci. Comput., 38 (2016),
  pp.~A3644--A3671.

\bibitem{Schumi}
{\sc L.~L. Schumaker}, {\em Spline functions: basic theory}, Cambridge
  Mathematical Library, Cambridge University Press, Cambridge, third~ed., 2007.

\bibitem{Sederberg_Zheng}
{\sc T.~Sederberg, J.~Zheng, A.~Bakenov, and A.~Nasri}, {\em T-splines and
  {T-NURCCS}s}, ACM Trans. Graph., 22 (2003), pp.~477--484.

\bibitem{GEOPDES-NEW}
{\sc R.~V{\'a}zquez}, {\em A new design for the implementation of isogeometric
  analysis in {O}ctave and {M}atlab: {G}eo{PDE}s 3.0}, Comput. Math. Appl., 72
  (2016), pp.~523 -- 554.

\bibitem{Vuong_giannelli_juttler_simeon}
{\sc A.-V. Vuong, C.~Giannelli, B.~J\"uttler, and B.~Simeon}, {\em A
  hierarchical approach to adaptive local refinement in isogeometric analysis},
  Comput. Methods Appl. Mech. Engrg., 200 (2011), pp.~3554--3567.

\bibitem{TT-splines}
{\sc X.~Wei, Y.~Zhang, L.~Liu, and T.~J. Hughes}, {\em Truncated {T}-splines:
  {F}undamentals and methods}, Comput. Methods Appl. Mech. Engrg.,  (2016),
  pp.~--.
\newblock In press.

\bibitem{Xu1992}
{\sc J.~Xu}, {\em Iterative methods by space decomposition and subspace
  correction}, SIAM Rev., 34 (1992), pp.~581--613.

\bibitem{Xu_MGintro}
\leavevmode\vrule height 2pt depth -1.6pt width 23pt, {\em An introduction to
  multigrid convergence theory}, in Iterative methods in scientific computing
  ({H}ong {K}ong, 1995), Springer, Singapore, 1997, pp.~169--241.

\bibitem{XCN}
{\sc J.~Xu, L.~Chen, and R.~H. Nochetto}, {\em Optimal multilevel methods for
  {$H({\rm grad})$}, {$H({\rm curl})$}, and {$H({\rm div})$} systems on graded
  and unstructured grids}, in Multiscale, nonlinear and adaptive approximation,
  Springer, Berlin, 2009, pp.~599--659.

\end{thebibliography}

\end{document}